\documentclass[12pt,a4paper]{article}
\usepackage[utf8x]{inputenc}
\usepackage{amsmath}
\usepackage{amsfonts}
\usepackage{amssymb}
\usepackage{amsthm}
\usepackage{mathrsfs}
\usepackage{fancyhdr}
\usepackage{graphicx}
\usepackage[usenames,x11names]{xcolor}
\usepackage{arabtex}
\usepackage{framed}
\usepackage[top=2cm,bottom=2cm,margin=2cm]{geometry}
\usepackage{cancel}
\usepackage{array}

\usepackage{hyperref}
\usepackage[numbers,sort&compress]{natbib}  

\title{Expansion of a bivariate symmetric mean in the neighborhood of the first bisector}
\author{\sc Bakir FARHI \\
National Higher School of Mathematics \\
P.O.Box 75, Mahelma 16093, Sidi Abdellah (Algiers) \\
Algeria \\[1mm]
\href{mailto:bakir.farhi@nhsm.edu.dz}{\tt bakir.farhi@nhsm.edu.dz} \\[1mm]
\url{http://farhi.bakir.free.fr/}
}
\date{}

\let\up=\textsuperscript
\let\epsilon=\varepsilon

\def\R{{\mathbb R}}

\def\N{{\mathbb N}}

\def\L{\mathscr{L}}
\def\I{\mathscr{I}}

\def\id{\mathrm{id}}

\def\AGM{\mathrm{AGM}}
\def\im{\mathrm{Im}}

\def\idem{\leavevmode\hbox to 10.6mm{\vrule height .63ex depth -.59ex
    width 10mm\hfill}}

\newcommand{\vabs}[1]{\left\vert{#1}\right\vert}

\theoremstyle{plain}
\numberwithin{equation}{section}
\newtheorem{thm}{Theorem}[section]

\newtheorem{lemma}[thm]{Lemma}
\newtheorem{prop}[thm]{Proposition}
\newtheorem{coll}[thm]{Corollary}

\theoremstyle{definition}
\newtheorem{defi}[thm]{Definition}
\newtheorem{defis}[thm]{Definitions}
\newtheorem{nota}[thm]{Notation}

\theoremstyle{remark}
\newtheorem{rmk}[thm]{Remark}
\newtheorem{rmks}[thm]{Remarks}
\newtheorem{expl}[thm]{Example}
\newtheorem{expls}[thm]{Examples}

\begin{document}
\maketitle

\begin{abstract}
In this paper, we investigate the behavior of a bivariate mean $M$ near the first bisector by establishing, in several significant cases, an important expansion of $M$ derived from the Taylor expansion of a single-variable function. These expansions are made explicit for a number of classical means. This motivates the introduction of the concept of the \emph{characteristic function} $Q_M$ of a mean $M$, defined as the second partial derivative of $M$ with respect to its first variable, evaluated along the diagonal. The function $Q_M$ measures the proximity of $M$ to the arithmetic mean near the first bisector and provides a univariate analytic framework for comparing and classifying means. We prove that inequalities between characteristic functions yield local inequalities between the corresponding means, and that in the case of homogeneous means, such inequalities hold globally. 

We also examine several important classes of means, both classical and novel, including: \emph{normal means}, \emph{additive means}, \emph{integral means of the first kind}, \emph{integral means of the second kind}, \emph{weighted integral means of the first kind}, and \emph{weighted integral means of the second kind}. For each class, we determine the specific form taken by the characteristic functions $Q_M$ of the means $M$ it contains, and we then study the injectivity and the surjectivity of the mapping $M \mapsto Q_M$ within the class. We also use characteristic functions to investigate intersections between certain classes of means, highlighting one of the key strengths of this concept. Finally, we introduce and study, for a given mean $M$, the class of \emph{$M$-means}, and show, in particular, that the arithmetic-geometric mean $\AGM$ is an $M$-mean for a specific weighted integral mean of the first kind $M$.   
\end{abstract}


\tableofcontents

\pagebreak

\section{Introduction and Notation}\label{sec1}
Throughout this paper, a \emph{mean} (or a \emph{bivariate mean}) refers to a function $M : (0 , + \infty)^2 \rightarrow (0 , + \infty)$ that is of class $C^{\infty}$ and satisfies the following two properties:
\begin{enumerate}
\item[(i)] $\forall x , y > 0$: 
$$
\min(x , y) \leq M(x , y) \leq \max(x , y) .
$$
\item[(ii)] $\forall x , y > 0$: 
$$
M(x , y) = M(y , x) ,
$$
that is, $M$ is symmetric.
\end{enumerate}
Among the most commonly used means, we mention the following (this list is not exhaustive):
\begin{itemize}
\item The \emph{Arithmetic mean}, denoted by $A$, and defined by:
$$
A(x , y) := \frac{x + y}{2} ~~~~~~~~~~ (\forall x , y > 0) .
$$
\item The \emph{Geometric mean}, denoted by $G$, and defined by:
$$
G(x , y) := \sqrt{x y} ~~~~~~~~~~ (\forall x , y > 0) .
$$
\item The \emph{Harmonic mean}, denoted by $H$, and defined by:
$$
H(x , y) := \frac{2 x y}{x + y} ~~~~~~~~~~ (\forall x , y > 0) .
$$
\item The \emph{$r$\up{th} Power mean} ($r \in \R^*$), denoted by $\pi_r$, and defined by:
$$
\pi_r(x , y) := \left(\frac{x^r + y^r}{2}\right)^{1/r} ~~~~~~~~~~ (\forall x , y > 0) .
$$
(See \cite{bul}). For $r = 2$, $\pi_r$ is called the \emph{quadratic mean} and simply denoted by $Q$. Thus, for all $x , y > 0$:
$$
Q(x , y) := \sqrt{\frac{x^2 + y^2}{2}} .
$$
Moreover, it is straightforward to verify that $\lim_{r \rightarrow 0} \pi_r = G$. Hence, it is customary to define (by continuity) $\pi_0 := G$. Additionally, note the following special cases:
$$
\pi_1 = A , \quad \pi_{-1} = H , \quad \pi_{1/2} = \frac{A + G}{2} .
$$
\item The \emph{Logarithmique mean}, denoted by $L$, and defined by:
$$
L(x , y) := \begin{cases}
x & \text{if } x = y \\
\dfrac{x - y}{\log{x} - \log{y}} & \text{otherwise}
\end{cases} \quad (\forall x , y > 0) .
$$
Several elegant formulas can also serve as definitions (or characterizations) of the logarithmic mean (see \cite{bul}):
\begin{align*}
L(x , y) & = \left(\int_{0}^{+ \infty} \dfrac{d t}{(t + x) (t + y)}\right)^{-1} , \\[1mm]
L(x , y) & = \int_{0}^{1} x^t y^{1 - t} \, d t , \\[1mm]
L(x , y) & = \prod_{n = 1}^{+ \infty} A\left(x^{1/2^n} , y^{1/2^n}\right)
\end{align*}
(for all $x , y > 0$).
\item The \emph{Exponential mean}, denoted by $E$, and defined by:
$$
E(x , y) := \begin{cases}
x & \text{if } x = y \\
\log\left(\dfrac{e^x - e^y}{x - y}\right) & \text{otherwise}
\end{cases} \quad (\forall x , y > 0) .
$$
\item The \emph{Lehmer mean}, denoted by $\L_t$ ($t \in \R$), and defined by:
$$
\L_t(x , y) := \frac{x^t + y^t}{x^{t - 1} + y^{t - 1}} \quad (\forall x , y > 0) .
$$
Note the following special cases:
$$
\L_1 = A , \quad \L_0 = H , \quad \L_{\frac{1}{2}} = G .
$$
(See \cite[Chapter 8]{bor} or \cite{bul}).
\item The \emph{Gini mean}, denoted by $\eta_{a , b}$ ($a , b \in \R$, $a \neq b$), and defined by:
$$
\eta_{a , b}(x , y) := \left(\frac{x^a + y^a}{x^b + y^b}\right)^{\frac{1}{a - b}} ~~~~~~~~~~ (\forall x , y > 0) .
$$
(See \cite[Chapter 8]{bor} or \cite{bul}). For $a \in \R$, we also define $\eta_{a , a} := \lim_{b \rightarrow a} \eta_{a , b}$ which simplifies to:
\begin{equation}\label{eq-intro-2}
\eta_{a , a}(x , y) = \exp\left(\frac{(\log x) x^a + (\log y) y^a}{x^a + y^a}\right) \quad (\forall x , y > 0) .
\end{equation}
Notice that the Gini means generalize both the power means and the Lehmer means. Specifically, for all $r , t \in \R$:
\begin{equation}\label{eq-intro-1}
\pi_r = \eta_{r , 0} \quad \text{and} \quad \L_t = \eta_{t , t - 1} .
\end{equation}
Furthermore, the Gini means can be expressed in terms of the Lehmer means as:
$$
\eta_{a , b} = \left(\L_{\frac{a}{a - b}}\left(x^{a - b} , y^{a - b}\right)\right)^{\frac{1}{a - b}}
$$
($\forall a , b \in \R$, with $a \neq b$, $\forall x , y > 0$).
\item The \emph{Stolarsky mean}, denoted by $S_p$ ($p \in \R \setminus \{0 , 1\}$), and defined by:
$$
S_p(x , y) := \begin{cases}
x & \text{if } x = y \\
\left(\dfrac{x^p - y^p}{p (x - y)}\right)^{\frac{1}{p - 1}} & \text{otherwise}
\end{cases} \quad (\forall x , y > 0) .
$$
It is easy to verify that $S_2 = A$, $S_{-1} = G$, and $\lim_{p \rightarrow 0} S_p = L$. Additionally, the limit case $\lim_{p \rightarrow 1} S_p$ gives the so-called \emph{Identric mean}, denoted by $I$, and is given by:
$$
I(x , y) = \begin{cases}
x & \text{if } x = y \\
e^{-1} \left(\dfrac{x^x}{y^y}\right)^{\frac{1}{x - y}} & \text{otherwise}
\end{cases} \quad (\forall x  , y > 0) .
$$
(See \cite[Chapter 8]{bor} or \cite{bul}). Hence, it is conventional to define (by continuity) $S_0 := L$ and $S_1 := I$.
\item The \emph{Gauss arithmetic-geometric mean}, denoted by $\AGM$, and defined as follows: Given two positive real numbers $x$ and $y$, the value $\AGM(x , y)$ is the common limit of the two sequences ${(x_n)}_{n \in \N_0}$ and ${(y_n)}_{n \in \N_0}$, which are recursively defined by
$$
\left\{\begin{array}{l}
x_0 = x , y_0 = y , \\[2mm]
x_{n + 1} = \dfrac{x_n + y_n}{2} \quad (\forall n \in \N_0) , \\[2mm]
y_{n + 1} = \sqrt{x_n y_n} \quad (\forall n \in \N_0) .
\end{array}
\right.
$$
Gauss discovered and established that the arithmetic-geometric mean is deeply related to elliptic integrals (see e.g., \cite{bor}). He explicitly proved the remarkable formula:
\begin{equation}\label{eq7}
\AGM(x , y) = \frac{\pi}{2} \left(\int_{0}^{+ \infty}\!\!\!\dfrac{d t}{\sqrt{(t^2 + x^2) (t^2 + y^2)}}\right)^{\!\!-1} \!\!\!= \frac{\pi}{2} \left(\int_{0}^{\frac{\pi}{2}}\!\!\dfrac{d \theta}{\sqrt{x^2 \cos^2 \theta + y^2 \sin^2 \theta}}\right)^{\!\!-1}
\end{equation}
(valid for all $x , y > 0$).
\end{itemize}

There are several important properties that a mean $M$ can satisfies. One of the most important is \emph{homogeneity}, which means that $M$ scales proportionality:
$$
M(\lambda x , \lambda y) = \lambda M(x , y) \quad (\forall \lambda , x , y > 0) .
$$
Next, a mean $M$ is said to be \emph{additively homogeneous} if it satisfies the following condition:
$$
M(x + \lambda , y + \lambda) = M(x , y) + \lambda \quad (\forall \lambda , x , y > 0) .
$$
In this case, $M$ can be naturally extended to a function on $\R^2$. Another property is ``monotonicity''. A mean $M$ is said to be \emph{isotone} (resp. \emph{strictly isotone}) if, for any fixed $y > 0$, the function $x \mapsto M(x , y)$ is nondecreasing (resp. increasing) on $(0 , + \infty)$ \cite[Chapter 8]{bor}. Among the previously mentioned means, $A$ and $E$ are additively homogeneous, while all except $E$ are homogeneous. We will see later that the arithmetic mean $A$ is the only mean that is both homogeneous and additively homogeneous. Furthermore, the Lehmer mean $\L_t$ (for $0 \leq t \leq 1$) and the arithmetic-geometric mean $\AGM$ are strictly isotone, whereas the Lehmer mean $\L_t$ fails to be isotone for $t > 1$. Strictness is another key property that plays a crucial role in the theory of means. A mean $M$ is said to be \emph{strict} (see \cite[Chapter 8]{bor}) if it satisfies the following property:
$$
\forall x , y > 0 :~ M(x , y) = x \Longrightarrow x = y .
$$
We can verify (although this is not always immediate) that all the means considered above are strict. The strictness property is so important that it is considered in \cite{far} as an integral part of the very definition of a mean. 

Given a mean $M$ and an infinitely differentiable function $f : (0 , + \infty) \rightarrow (0 , + \infty)$, we define the two variable functions $f \circ M$ and $M \circ f$ by:
\begin{equation*}
\begin{split}
\left(f \circ M\right)(x , y) & := f\left(M(x , y)\right) , \\
\left(M \circ f\right)(x , y) & := M\left(f(x) , f(y)\right) 
\end{split} \quad (\forall x , y > 0) .
\end{equation*}
When $f$ happens to be strictly monotonic, the composition $f^{-1} \circ M \circ f$ defines a new mean. This provides one of the most common ways for generating new means from given ones.

The theory of means has long been a subject of interest to both classical and contemporary authors. Two major problems have particularly drawn attention: the comparison between means, and the determination of the intersection between certain classes of means. To address these questions, researchers --- Starting from Lehmer \cite{leh} --- have employed various types of expansions of means. For homogeneous means $M$, Lehmer \cite{leh} and later Alzer and Ruscheweyh \cite{alz} considered the power series expansion of the single-variable function $M(1 , 1 - t)$ for $0 \leq t < 1$. More recently, in a series of papers \cite{bu-el,el1,elv1,elv2,elv3}, Buri\`c, Elezovi\`c, and Mihokovi\`c developed the concept of asymptotic expansions of homogeneous means and obtained significant results regarding the comparison of bivariate means and the intersection of certain classes of homogeneous means. Their approach is based on representing the function
$$
t \mapsto M(t + x , t + y) = t \, M\left(1 + \frac{x}{t} , 1 + \frac{y}{t}\right)
$$
as $t \to + \infty$ in the form
$$
t + c_1(x , y) + \dfrac{c_2(x , y)}{t} + \dfrac{c_3(x , y)}{t^2} + \dots ,
$$
which arises from the Taylor expansion of $M$ at the point $(1 , 1)$. For instance, using this method, Mihokovi\`c \cite{mih} provided a partial solution to an important problem posed by the author in \cite{far}, concerning the coincidence of two symmetries of different nature within the set of means.

The aim of the present paper is to introduce a new type of expansion, valid for means not necessarily homogeneous, which turns out to be highly effective in addressing both aforementioned problems. This is the expansion of a bivariate mean near the first bisector. It naturally leads to the notion of the \emph{characteristic function} $Q_M$ of a mean $M$, defined as the second partial derivative of $M$ with respect to its first variable, evaluated along the diagonal. We demonstrate that this concept provides a powerful univariate analytic tool for comparing means near the first bisector and, more remarkably, for comparing homogeneous means globally. Furthermore, we examine several classes of means, some well-known and other news, including: \emph{normal means}, \emph{additive means}, \emph{integral means of the first kind}, \emph{integral means of the second kind}, \emph{weighted integral means of the first kind}, and \emph{weighted integral means of the second kind}. For each of these classes, we determine the specific form taken by the characteristic functions $Q_M$ of the means $M$ it contains, and we then study their \emph{minimality} and \emph{completeness}, that is, the injectivity and surjectivity of the mapping $M \mapsto Q_M$ within the class. We also use characteristic functions to investigate intersections between certain classes of means, illustrating one of the conceptual strengths of this tool. As a final contribution, we introduce the notion of \emph{$M$-means} associated with a given mean $M$, explore some of their structural properties, and establish that the arithmetic-geometric mean $\AGM$ arises as an $M$-mean for a suitable weighted integral mean of the first kind $M$.

\section{Expansion of a mean near the first bisector}

\subsection{The general case}

Let $M$ be a mean. We associate with $M$ the family of real-valued functions $\varphi_x$ ($x \in (0 , + \infty)$), defined by:
\begin{equation}\label{eq3}
\begin{array}{rcl}
\varphi_x :~ (- 2 x , 2 x) & \longrightarrow & (0 , + \infty) \\
t & \longmapsto & \varphi_x(t) := M\left(x + \frac{t}{2} , x - \frac{t}{2}\right)
\end{array} .
\end{equation}
Then, using the properties of $M$ as a mean, we have for all $x \in (0 , + \infty)$ and all $t \in (- 2 x , 2 x)$:
\begin{align*}
\varphi_x(0) & = M(x , x) = x \\[-6mm]
\intertext{and} \\[-12mm]
\varphi_x(-t) & = M\left(x - \frac{t}{2} , x + \frac{t}{2}\right) = M\left(x + \frac{t}{2} , x - \frac{t}{2}\right) = \varphi_x(t) .
\end{align*}
In particular, the last property shows that the functions $\varphi_x$ are all even. Assuming that the functions $\varphi_x$ admit a Taylor series expansion around $0$, we can write for $x \in (0 , + \infty)$ and $t$ in a neighborhood of $0$:
$$
\varphi_x(t) = x + \frac{\varphi_x''(0)}{2!} t^2 + \frac{\varphi_x^{(4)}(0)}{4!} t^4 + \dots .
$$
Setting for simplicity
\begin{equation}\label{eq1}
f_n(x) := \frac{\varphi_x^{(2 n)}(0)}{(2 n)!} \quad (\forall n \in \N) ,
\end{equation}
we obtain
$$
\varphi_x(t) = x + f_1(x) t^2 + f_2(x) t^4 + \dots .
$$
Finally, if $x , y \in (0 , + \infty)$ are sufficiently close to each other, substituting $x$ by $A = A(x , y) = \frac{x + y}{2}$ and $t$ by $x - y$, we get
$$
M(x , y) = A + f_1(A) (x - y)^2 + f_2(A) (x - y)^4 + \dots .
$$
Thus, we have proven the following theorem:
\begin{thm}\label{t1}
Let $M$ be a mean such that all associated functions $\varphi_x$ ($x > 0$), given in \eqref{eq3}, admit a Taylor series expansion around $0$. Then, there exists a sequence of functions $(f_n)_{n \in \N}$, defined on $(0 , + \infty)$, such that for $(x , y)$ in a neighborhood of the first bisector:
\begin{equation}\label{eq2}
M(x , y) = A + f_1(A) (x - y)^2 + f_2(A) (x - y)^4 + \dots .
\end{equation}
\end{thm}

\noindent\textbf{Terminology:} In the context of Theorem \ref{t1}, Equation \eqref{eq2} is called \emph{the expansion of the mean $M$ near the first bisector} (or, for brevity, simply \emph{the expansion of $M$}).

\begin{rmks}~
\begin{enumerate}
\item By the uniqueness of the Taylor series expansion around $0$ for a real-valued function of one variable, it follows that if a mean admits an expansion near the first bisector, then this expansion is necessarily unique.
\item Equation \eqref{eq2} in Theorem \ref{t1} shows that the affine approximation of any (regular) mean near the first bisector is the arithmetic mean.
\end{enumerate}
\end{rmks}

\subsection{The case of a homogeneous mean}

If a mean $M$ is homogeneous, then using the above considerations, for all $x \in (0 , + \infty)$ and all $t \in (- 2 x , 2 x)$, we obtain:
$$
\varphi_x(t) := M\left(x + \frac{t}{2} , x - \frac{t}{2}\right) = x \, M\left(1 + \frac{t}{2 x} , 1 - \frac{t}{2 x}\right) = x \, \varphi_1\left(\frac{t}{x}\right) .
$$
Thus, the Taylor series expansions of the functions $\varphi_x$ around $0$ (if they exist) are all derived from the Taylor series expansion of $\varphi_1$ around $0$. Consequently, the expansion of $M$ in the neighborhood of the first bisector (if it exists) can be directly obtained from the Taylor series expansion of $\varphi_1$, as specified in the following corollary:
\begin{coll}\label{c1}
Let $M$ be a homogeneous mean. Suppose that the function $\varphi_1$ defined by:
$$
\begin{array}{rcl}
\varphi_1 :~ (-2 , 2) & \longrightarrow & (0 , + \infty) \\
t & \longmapsto & \varphi_1(t) := M\left(1 + \frac{t}{2} , 1 - \frac{t}{2}\right)
\end{array}
$$
admits a Taylor series expansion around $0$:
$$
\varphi_1(t) = 1 + c_1 t^2 + c_2 t^4 + \dots ,
$$
where $c_1 , c_2 , \dots$ are real coefficients. Then, for all sufficiently close positive real numbers $x$ and $y$, we have:
$$
M(x , y) = A + \frac{c_1}{A} (x - y)^2 + \frac{c_2}{A^3} (x - y)^4 + \frac{c_3}{A^5} (x - y)^6 + \dots .
$$
\end{coll}

\begin{proof}
Setting $t = \frac{x - y}{A}$ in the Taylor series expansion of $\varphi_1$ and using the homogeneity of $M$ yields the result. 
\end{proof}

\begin{rmk}
Conversely to Corollary \ref{c1}, we can immediately verify that if the expansion of a mean $M$ has the form
$$
M(x , y) = A + \frac{c_1}{A} (x - y)^2 + \frac{c_2}{A^3} (x - y)^4 + \frac{c_3}{A^5} (x - y)^6 + \dots
$$
(where $c_1 , c_2 , c_3 , \dots$ are real numbers) then $M$ is necessarily homogeneous. So, an expansion of the form  
$$
A + \sum_{n = 1}^{+ \infty} \frac{c_n}{A^{2 n - 1}} (x - y)^{2 n}, \quad (c_n \in \R, \forall n \in \N),
$$  
for a mean \( M \) characterizes its homogeneity.
\end{rmk}

\subsection{The case of an additively homogeneous mean}

If a mean $M$ is additively homogeneous, then it naturally extends to $\R^2$. Consequently, the functions $\varphi_x$ ($x > 0$) associated to $M$ and defined in \eqref{eq3}, simplify to
$$
\varphi_x(t) := M\left(x + \frac{t}{2} , x - \frac{t}{2}\right) = x - \frac{t}{2} + M(t , 0) \quad (\forall t \in \R) .
$$
Thus, if they exist, the Taylor series expansions of the functions $\varphi_x$ around $0$ are directly obtained from the Taylor series expansion of the function $t \mapsto M(t , 0)$ around $0$. Consequently, the expansion of $M$ near the first bisector (if it exists) follows directly from the Taylor series expansion of $t \mapsto M(t , 0)$, which necessarily takes the form
$$
\frac{t}{2} + \sum_{n = 1}^{+ \infty} c_n t^{2 n}
$$
(where the $c_n$'s are real numbers), because the additive homogeneity and symmetry of $M$ ensure that the function $M(t , 0) - \frac{t}{2} = M(\frac{t}{2} , - \frac{t}{2})$ is even. This is specified in the following corollary:
\begin{coll}\label{c2}
Let $M$ be an additively homogeneous mean. Suppose that the function $t \mapsto M(t , 0)$ admits a Taylor series expansion around $0$, given by
$$
M(t , 0) = \frac{1}{2} t + \sum_{n = 1}^{+ \infty} c_n t^{2 n} ,
$$
where $c_1 , c_2 , \dots$ are real coefficients. Then, for all sufficiently close positive real numbers $x$ and $y$, we have:
$$
M(x , y) = A + \sum_{n = 1}^{+ \infty} c_n (x - y)^{2 n} .
$$
\end{coll} 

\begin{proof}
For all sufficiently close positive real numbers $x , y$, we have
$$
M(x , y) = y + M(x - y , 0) = y + \frac{1}{2} (x - y) + \sum_{n = 1}^{+ \infty} c_n (x - y)^{2 n} = A(x , y) + \sum_{n = 1}^{+ \infty} c_n (x - y)^{2 n} ,
$$
as required.
\end{proof}

\begin{rmk}
Notably, the functions $f_1 , f_2 , \dots$ appearing in the expansion \eqref{eq2} of a general mean $M$ become constant if $M$ is additively homogeneous. Moreover, the converse of this property is also true by direct verification. Therefore, an expansion of the form $A + \sum_{n = 1}^{+ \infty} c_n (x - y)^{2 n}$ ($c_n \in \R$ for all $n \in \N$) for a mean characterizes its additive homogeneity. 
\end{rmk}

\subsection{Application to some usual means}

Applying Corollaries \ref{c1} and \ref{c2}, we derive the explicit expansions of some common homogeneous and additively homogeneous means. The following corollary provides examples for the three means $H , G$, and $E$.

\begin{coll}\label{c3}
The expansions of the means $H , G$, and $E$ are given respectively by:
\begin{align}
H(x , y) & = A - \frac{1}{4 A} (x - y)^2 , \label{eq4} \\[1mm]
G(x , y) & = A - \sum_{n = 1}^{+ \infty} \frac{\binom{2 n}{n}}{16^n (2 n - 1) A^{2 n - 1}} (x - y)^{2 n} , \label{eq5} \\[1mm]
E(x , y) & = A + \sum_{n = 1}^{+ \infty} \frac{B_{2 n}}{(2 n) \cdot (2 n)!} (x - y)^{2 n} , \label{eq6}
\end{align}
where $B_k$ are the Bernoulli numbers, defined via their exponential generating function:
$$
\frac{t}{e^t - 1} = \sum_{k = 0}^{+ \infty} B_k \frac{t^k}{k!} = 1 - \frac{t}{2} + \sum_{k = 0}^{+ \infty} B_{2 k} \frac{t^{2 k}}{(2 k)!} .
$$
\end{coll}

\begin{proof}
For simplicity and clarity, let us denote by $h , g$ the functions $\varphi_1$ in Corollary \ref{c1} associated respectively to the means $H$ and $G$, and let $e$ denote the function $t \mapsto E(t , 0)$. First, for all $t \in (-2 , 2)$, we have:
$$
h(t) := H\left(1 + \frac{t}{2} , 1 - \frac{t}{2}\right) = \frac{2 \left(1 + \frac{t}{2}\right) \left(1 - \frac{t}{2}\right)}{\left(1 + \frac{t}{2}\right) + \left(1 - \frac{t}{2}\right)} = 1 - \frac{t^2}{4} .
$$
Applying Corollary \ref{c1} immediately gives the expansion of $H$ in \eqref{eq4}.

Next, for all $t \in (-2 , 2)$, we have:
$$
g(t) := G\left(1 + \frac{t}{2} , 1 - \frac{t}{2}\right) = \sqrt{\left(1 + \frac{t}{2}\right) \left(1 - \frac{t}{2}\right)} = \left(1 - \frac{t^2}{4}\right)^{1/2} .
$$
Using the generalized binomial theorem:
$$
\left(1 - \frac{t^2}{4}\right)^{1/2} = \sum_{n = 0}^{+ \infty} \binom{1/2}{n} \left(- \frac{t^2}{4}\right)^n ,
$$
and noting that for all $n \in \N$,
$$
\binom{1/2}{n} = (-1)^{n - 1} \frac{\binom{2 n}{n}}{4^n (2 n - 1)} ,
$$
it follows that
$$
g(t) = 1 - \sum_{n = 1}^{+ \infty} \frac{\binom{2 n}{n}}{16^n (2 n -1)} t^{2 n} .
$$
Applying Corollary \ref{c1} yields the expansion of $G$ in \eqref{eq5}.

Finally, for all $t \in \R$, we have
$$
e(t) := E(t , 0) = \log\left(\frac{e^t - 1}{t}\right) .
$$
To apply Corollary \ref{c2}, we expand $t \mapsto \log(\frac{e^t - 1}{t})$ around $0$. We have
\begin{align*}
t \frac{d}{dt}\!\!\left(\log\left(\frac{e^t - 1}{t}\right)\right) & = t \left(\dfrac{\frac{d}{d t}(e^t - 1)}{e^t - 1} - \dfrac{\frac{d}{d t}(t)}{t}\right) \\
& = t \left(\dfrac{e^t}{e^t - 1} - \frac{1}{t}\right) \\
& = t - 1 + \dfrac{t}{e^t - 1} \\
& = \frac{t}{2} + \sum_{n = 1}^{+ \infty} B_{2 n} \frac{t^{2 n}}{(2 n)!}
\end{align*}
(for $t$ sufficiently close to $0$). Dividing by $t$ and integrating from $0$ to $t$ gives
$$
\log\left(\dfrac{e^t - 1}{t}\right) = \frac{t}{2} + \sum_{n = 1}^{+ \infty} \frac{B_{2 n}}{(2 n) \cdot (2 n)!} t^{2 n} .
$$
Applying Corollary \ref{c2} yields the expansion of the mean $E$ in \eqref{eq6} and completes this proof.
\end{proof}

\subsection{Means that are both homogeneous and additively homogeneous}

Comparing Corollaries \ref{c1} and \ref{c2} and using the uniqueness of the expansion of a mean near the first bisector, we immediately obtain the following result.

\begin{coll}
The unique mean $M$ that admits an expansion near the first bisector and is both homogeneous and additively homogeneous is the arithmetic mean $A$. \hfill $\square$
\end{coll}

In fact, this corollary can be generalized to a broader class of real-valued functions on $(0 , + \infty)^2$. More precisely, we have the following result:

\begin{prop}\label{p1}
Let $f$ be a real-valued function defined on $(0 , + \infty)^2$, satisfying the following conditions:
\begin{enumerate}
\item[(i)] $f(1 , 1) = 1$.
\item[(ii)] $f$ admits partial derivatives at $(1 , 1)$.
\item[(iii)] $f$ is symmetric, i.e., $f(x , y) = f(y , x)$ for all $x , y > 0$.
\item[(iv)] $f$ is homogeneous: $f(\lambda x , \lambda y) = \lambda f(x , y)$ for all $\lambda , x , y > 0$.
\item[(v)] $f$ is additively homogeneous: $f(x + \lambda , y + \lambda) = f(x , y) + \lambda$ for all $\lambda , x , y > 0$.
\end{enumerate}
Then $f$ must necessarily be the arithmetic mean, i.e., $f = A$.
\end{prop}

\begin{proof}
Using Conditions (iv) and (v), we have for all $x , y , h > 0$:
$$
f(x , y) = \dfrac{f(1 + h x , 1 + h y) - 1}{h} .
$$
Taking the limit as $h \to 0^+$ and using Condition (ii), we obtain for all $x , y > 0$:
$$
f(x , y) = \frac{\partial f}{\partial x}(1 , 1) \cdot x + \frac{\partial f}{\partial y}(1 , 1) \cdot y .
$$
Since $f$ is symmetric by (iii), we have $\frac{\partial f}{\partial x}(1 , 1) = \frac{\partial f}{\partial y}(1 , 1)$. So, setting $a := \frac{\partial f}{\partial x}(1 , 1) = \frac{\partial f}{\partial y}(1 , 1)$, we deduce that for all $x , y > 0$:
$$
f(x , y) = a (x + y) .
$$
Finally, setting $x = y = 1$ in this formula and using (i), we obtain $a = \frac{1}{2}$, which leads to $f(x , y) = \frac{x + y}{2} = A(x , y)$ for all $x , y > 0$. This completes the proof.
\end{proof}

\begin{rmk}
Proposition \ref{p1} does not hold in general if Condition (ii) is omitted, as shown by the following counterexample:
$$
f(x , y) := \begin{cases}
\frac{1}{3} x + \frac{2}{3} y & \text{if } x \geq y , \\
\frac{2}{3} x + \frac{1}{3} y & \text{if } x < y
\end{cases} \quad (\text{for all } x , y > 0) .
$$
Here $f \neq A$ satisfies all the Conditions of Proposition \ref{p1}, except Condition (ii). We will see below that if we extend our domain to $\R^2$ instead of $(0 , + \infty)^2$, then Conditions (iii), (iv), and (v) (suitably adapted to $\R^2$) suffice to conclude that $f = A$.
\end{rmk}

\begin{prop}\label{p2}
Let $f$ be a real-valued function defined on $\R^2$, satisfying the following conditions:
\begin{enumerate}
\item[(i)] $f$ is symmetric, i.e., $f(x , y) = f(y , x)$ for all $x , y \in \R$.
\item[(ii)] $f$ is homogeneous: $f(\lambda x , \lambda y) = \lambda f(x , y)$ for all $\lambda , x , y \in \R$.
\item[(iii)] $f$ is additively homogeneous: $f(x + \lambda , y + \lambda) = f(x , y) + \lambda$ for all $\lambda , x , y \in \R$.
\end{enumerate}
Then $f$ must necessarily be the arithmetic mean, i.e., $f = A$.
\end{prop}

\begin{proof}
Using Conditions (ii) and (iii), we have for all $x , y \in \R$:
$$
f(x , y) = y + f(x - y , 0) = y + (x - y) f(1 , 0) .
$$
Setting $(x , y) = (0 , 1)$, we get
$$
f(0 , 1) = 1 - f(1 , 0) ,
$$
Since $f(0 , 1) = f(1 , 0)$ by (i), it follows that
$$
f(1 , 0) = \frac{1}{2} .
$$
Substituting this value of $f(1 , 0)$ into the previous formula of $f(x , y)$ gives
$$
f(x , y) = \frac{x + y}{2} = A(x , y) \quad (\forall x , y \in \R) ,
$$
as required.
\end{proof}

\subsection{The characteristic functions of a mean}

\subsubsection{The concept of the characteristic functions for a mean}

Given a mean $M$ that admits an expansion near the first bisector, we are now interested in expressing the functions $f_n$ in Formula \eqref{eq2} in terms of the higher-order partial derivatives of $M$ at points on the first bisector. Before doing so, we first establish that the functions $\frac{\partial^n M}{{\partial x}^i {\partial y}^j}(x , x)$ ($n , i , j \in \N_0$, $i + j = n$) can all be expressed as linear combinations (with rational coefficients independent of $M$) of the basic functions $\frac{\partial^{2 k} M}{{\partial x}^{2 k}}(x , x)$ ($k \in \N_0$) and their higher-order derivatives. We refer to these basic functions as the \emph{characteristic functions} of the mean $M$.

\begin{defi}\label{defi1}
Let $M$ be a mean. 
\begin{itemize}
\item For $k \in \N_0$, the real function $Q_M^{[k]}$ defined on $(0 , + \infty)$ by
$$
Q_M^{[k]}(x) := \dfrac{\partial^{2 k} M}{{\partial x}^{2 k}}(x , x)
$$
is called the \emph{characteristic function of order $k$} of $M$.
\item If no ambiguity arises regarding the mean $M$, we may omit the index $M$ in $Q_M^{[k]}$ and simply write $Q^{[k]}$.
\item Clearly, $Q_M^{[0]}(x) = x$.
\item The function $Q_M^{[1]}$ is denoted simply by $Q_M$ and referred to as the \emph{characteristic function} of $M$. We will later see that comparing between two means $M_1$ and $M_2$ near the first bisector reduces to comparing their characteristic functions $Q_{M_1}$ and $Q_{M_2}$ (see \S \ref{subsubsec2-b}). 
\end{itemize}
\end{defi} 

We establish the following proposition:

\begin{prop}\label{p3}
Let $M$ be a mean. Then for all $n , k \in \N_0$ with $k \leq n$, there exist rational numbers $c_i$ ($0 \leq i \leq \lfloor\frac{n}{2}\rfloor$), independent of $M$, such that the following identity holds:
\begin{equation}\label{eq8}
\dfrac{\partial^n M}{{\partial x}^k {\partial y}^{n - k}}(x , x) = \sum_{0 \leq i \leq \frac{n}{2}} c_i \left(\left(\frac{d}{d x}\right)^{\!\! n - 2 i} Q_M^{[i]}\right)(x) .
\end{equation} 
\end{prop}

Given the intricate nature of the proof of this proposition, we prefer to explicitly verify its validity for the first few values of $n$

Let $M$ be a mean. Differentiating the identity $M(x , x) = x$ yields
$$
\frac{\partial M}{\partial x}(x , x) + \frac{\partial M}{\partial y}(x , x) = 1 .
$$
Since $M$ is symmetric, we have $\frac{\partial M}{\partial x}(x , x) = \frac{\partial M}{\partial y}(x , x)$, leading to
\begin{equation}\label{eq-qay-1}
\frac{\partial M}{\partial x}(x , x) = \frac{\partial M}{\partial y}(x , x) = \frac{1}{2} .
\end{equation}
Differentiating again the identity $\frac{\partial M}{\partial x}(x , x) = \frac{1}{2}$ (given in \eqref{eq-qay-1}), we obtain
$$
\frac{\partial^2 M}{{\partial x}^2}(x , x) + \frac{\partial^2 M}{\partial x \partial y}(x , x) = 0 .
$$
Using the symmetry of $M$, it follows that:
\begin{align}
\frac{\partial^2 M}{{\partial x}^2}(x , x) & = \frac{\partial^2 M}{{\partial y}^2}(x , x) = Q_M^{[1]}(x) , \label{eq-qay-2} \\
\frac{\partial^2 M}{\partial x \partial y}(x , x) & = - Q_M^{[1]}(x) . \label{eq-qay-3}
\end{align}
Then, starting from $\frac{\partial^2 M}{{\partial x}^2}(x , x) = Q_M^{[1]}(x)$ and differentiating, we obtain
\begin{equation}\label{eq-qay-4}
\frac{\partial^3 M}{{\partial x}^3}(x , x) + \frac{\partial^3 M}{{\partial x}^2 \partial y}(x , x) = \frac{d}{d x}\left(Q_M^{[1]}\right)(x) .
\end{equation}
Next, differentiating $\frac{\partial^2 M}{\partial x \partial y}(x , x) = - Q_M^{[1]}(x)$ (given in \eqref{eq-qay-3}), we get
$$
\frac{\partial^3 M}{{\partial x}^2 \partial y}(x , x) + \frac{\partial^3 M}{\partial x {\partial y}^2}(x , x) = - \frac{d}{d x}\left(Q_M^{[1]}\right)(x) .
$$
By the symmetry of $M$, it follows that:
\begin{equation}\label{eq-qay-5}
\frac{\partial^3 M}{{\partial x}^2 \partial y}(x , x) = \frac{\partial^3 M}{\partial x {\partial y}^2}(x , x) = - \frac{1}{2} \frac{d}{d x}\left(Q_M^{[1]}\right)(x) .
\end{equation}
Substituting this into \eqref{eq-qay-4}, and using the symmetry of $M$, we conclude that:
\begin{equation}\label{eq-qay-6}
\frac{\partial^3 M}{{\partial x}^3}(x , x) = \frac{\partial^3 M}{{\partial y}^3}(x , x) = \frac{3}{2} \frac{d}{d x}\left(Q_M^{[1]}\right)(x) .
\end{equation}
Continuing in the same manner, differentiating $\frac{\partial^3 M}{{\partial x}^3}(x , x) = \frac{3}{2} \frac{d}{d x}\left(Q_M^{[1]}\right)(x)$, we obtain
$$
\frac{\partial^4 M}{{\partial x}^4}(x , x) + \frac{\partial^4 M}{{\partial x}^3 \partial y}(x , x) = \frac{3}{2} \left(\left(\frac{d}{d x}\right)^2 Q_M^{[1]}\right)(x) ,
$$
which gives (taking into account again the symmetry of $M$):
\begin{equation}\label{eq-qay-7}
\frac{\partial^4 M}{{\partial x}^3 \partial y}(x , x) = \frac{\partial^4 M}{\partial x {\partial y}^3}(x , x) = \frac{3}{2} \left(\left(\frac{d}{d x}\right)^2 Q_M^{[1]}\right)(x) - Q_M^{[2]}(x) .
\end{equation}
Similarly, differentiating $\frac{\partial^3 M}{{\partial x}^2 \partial y}(x , x) = - \frac{1}{2} \frac{d}{d x}\left(Q_M^{[1]}\right)(x)$ (given in \eqref{eq-qay-5}), we obtain
\begin{equation}\label{eq-qay-8}
\frac{\partial^4 M}{{\partial x}^3 \partial y}(x , x) + \frac{\partial^4 M}{{\partial x}^2 {\partial y}^2}(x , x) = - \frac{1}{2} \left(\left(\frac{d}{d x}\right)^2 Q_M^{[1]}\right)(x) .
\end{equation}
By substituting into \eqref{eq-qay-8} $\frac{\partial^4 M}{{\partial x}^3 \partial y}(x , x)$ by its expression given in \eqref{eq-qay-7}, we get
\begin{equation}\label{eq-qay-9}
\frac{\partial^4 M}{{\partial x}^2 {\partial y}^2}(x , x) = - 2 \left(\left(\frac{d}{d x}\right)^2 Q_M^{[1]}\right)(x) + Q_M^{[2]}(x) .
\end{equation}
Continuing this differentiation process, we also establish that
\begin{align}
\frac{\partial^5 M}{{\partial x}^5}(x , x) & = \frac{\partial^5 M}{{\partial y}^5}(x , x) = - \frac{5}{2} \left(\left(\frac{d}{d x}\right)^3 Q_M^{[1]}\right)(x) + \frac{5}{2} \frac{d}{d x} \left(Q_M^{[2]}\right)(x) , \label{eq-qay-10} \\
\frac{\partial^5 M}{{\partial x}^4 \partial y}(x , x) & = \frac{\partial^5 M}{\partial x {\partial y}^4}(x , x) = \frac{5}{2} \left(\left(\frac{d}{d x}\right)^3 Q_M^{[1]}\right)(x) - \frac{3}{2} \frac{d}{d x} \left(Q_M^{[2]}\right)(x) , \label{eq-qay-11} \\
\frac{\partial^5 M}{{\partial x}^3 {\partial y}^2}(x , x) & = \frac{\partial^5 M}{{\partial x}^2 {\partial y}^3}(x , x) = - \left(\left(\frac{d}{d x}\right)^3 Q_M^{[1]}\right)(x) + \frac{1}{2} \frac{d}{d x} \left(Q_M^{[2]}\right)(x) , \label{eq-qay-12} 
\end{align}
and so on. As we can see, Formulas \eqref{eq-qay-1}, \eqref{eq-qay-2}, \eqref{eq-qay-3}, \eqref{eq-qay-5}, \eqref{eq-qay-6}, \eqref{eq-qay-7}, \eqref{eq-qay-9}, \eqref{eq-qay-10}, \eqref{eq-qay-11}, and \eqref{eq-qay-12} all support the result of Proposition \ref{p3}.

We now proceed to prove Proposition \ref{p3} in its general case.

\begin{proof}[Proof of Proposition \ref{p3}]
We proceed by \textbf{nested induction}, where the outer induction is on $n \in \N_0$ and the inner induction is on $k \in \{0 , 1 , \dots , n\}$. \\[1mm]
\textbullet{} For small values of $n$, the formula of Proposition \ref{p3} has already been verified explicitly; namely, for $n \in \{0 , 1 , 2 , 3 , 4 , 5\}$. \\[1mm]
\textbullet{} Let $n \geq 6$ be an integer. Assume that Formula \eqref{eq8} of Proposition \ref{p3} holds for $(n - 1)$ and all $k \in \N_0$ with $k \leq n - 1$. We aim to show that it also holds for $n$ and all $k \in \N_0$ with $k \leq n$. To achieve this, we perform an inner induction on $k$, distinguishing two cases according to the parity of $n$. \\[1mm]
\textbf{Inductive step for the inner induction.} Assume that Formula \eqref{eq8} holds for some 
pair $(n , k)$ with $k \in \{0 , 1 , \dots , n - 1\}$ and prove that it must also hold for the pair $(n , k + 1)$. The outer induction hypothesis ensures the existence of rational numbers $a_i$ (independent of $M$) such that
$$
\dfrac{\partial^{n - 1} M}{{\partial x}^k {\partial y}^{n - 1 - k}}(x , x) = \sum_{0 \leq i \leq \frac{n - 1}{2}} a_i \left(\left(\frac{d}{d x}\right)^{\!\! n - 1 - 2 i} Q_M^{[i]}\right)(x) .
$$
Differentiating, we obtain
$$
\dfrac{\partial^n M}{{\partial x}^{k + 1} {\partial y}^{n - k - 1}}(x , x) + \dfrac{\partial^n M}{{\partial x}^k {\partial y}^{n - k}}(x , x) = \sum_{0 \leq i \leq \frac{n - 1}{2}} a_i \left(\left(\frac{d}{d x}\right)^{\!\! n - 2 i} Q_M^{[i]}\right)(x) .
$$
Since the inner induction hypothesis guarantees that $\frac{\partial^n M}{{\partial x}^k {\partial y}^{n - k}}(x , x)$ is a linear combination of the terms $\left(\frac{d}{d x}\right)^{n - 2 i} Q_M^{[i]}$ ($0 \leq i \leq \frac{n}{2}$) with rational coefficients independent of $M$, it follows that the same holds for $\frac{\partial^n M}{{\partial x}^{k + 1} {\partial y}^{n - k - 1}}(x , x)$, as required. This completes the inductive step for the inner induction.

Now, we handle the base case of the inner induction by distinguishing two cases: \\[1mm]
\textbf{Case 1: $n$ is even.} In this case, we have by the symmetry of $M$:
$$
\dfrac{\partial^n M}{{\partial y}^n}(x , x) = \dfrac{\partial^n M}{{\partial x}^n}(x , x) = Q_M^{[\frac{n}{2}]}(x) ,
$$ 
which directly verifies Formula \eqref{eq8} for the pair $(n , 0)$. By the inner inductive step, \eqref{eq8} holds for all the pairs $(n , k)$, with $k \in \{0 , 1 , \dots , n\}$. \\[1mm]
\textbf{Case 2: $n$ is odd.} Write $n = 2 m + 1$ for some $m \in \N$. By the outer induction hypothesis, there exist rational numbers $b_i$ (independent of $M$) such that
$$
\dfrac{\partial^{2 m} M}{{\partial x}^m {\partial y}^m}(x , x) = \sum_{0 \leq i \leq m} b_i \left(\left(\frac{d}{d x}\right)^{\!\! 2 m - 2 i} Q_M^{[i]}\right)(x) .
$$
Differentiating, we obtain
$$
\dfrac{\partial^{2 m + 1} M}{{\partial x}^{m + 1} {\partial y}^m}(x , x) + \dfrac{\partial^{2 m + 1} M}{{\partial x}^m {\partial y}^{m + 1}}(x , x) = \sum_{0 \leq i \leq m} b_i \left(\left(\frac{d}{d x}\right)^{\!\! 2 m + 1 - 2 i} Q_M^{[i]}\right)(x) .
$$
By symmetry of $M$, we have 
$$
\frac{\partial^{2 m + 1} M}{{\partial x}^m {\partial y}^{m + 1}}(x , x) = \frac{\partial^{2 m + 1} M}{{\partial x}^{m + 1} {\partial y}^m}(x , x) ,
$$
so that
$$
\dfrac{\partial^{2 m + 1} M}{{\partial x}^m {\partial y}^{m + 1}}(x , x) = \frac{1}{2} \sum_{0 \leq i \leq m} b_i \left(\left(\frac{d}{d x}\right)^{\!\! 2 m + 1 - 2 i} Q_M^{[i]}\right)(x) .
$$
This proves Formula \eqref{eq8} for the pair $(n , m) = (n , \frac{n - 1}{2})$. By the inner inductive step, Formula \eqref{eq8} must hold for all the pairs $(n , k)$, with $\frac{n - 1}{2} \leq k \leq n$. Finally, by symmetry of $M$, we conclude that it holds for all the pairs $(n , k)$, with $k \in \{0 , 1 , \dots , n\}$, since
$$
\dfrac{\partial^n M}{{\partial x}^{n - k} {\partial y}^k}(x , x) = \dfrac{\partial^n M}{{\partial x}^k {\partial y}^{n - k}}(x , x) \quad (\forall k \in \{0 , 1 , \dots , n\}) .
$$
This completes the nested induction and the proof.
\end{proof}

\begin{rmk}
The explicit closed-form expression of the rational coefficients appearing in Formula \eqref{eq8} of Proposition \ref{p3} is not straightforward and remains an open question. However, these coefficients can be computed recursively, as implicitly described in the proof of Proposition \ref{p3}.
\end{rmk}

Below, we present an elegant polynomial formulation capturing the relationships between partial derivatives of various orders of a given mean $M$ at first bisector points. This is given by the following proposition:

\begin{prop}
Let $M$ be a mean. For all $n \in \N_0$, define
$$
F_n(x ,T) := \sum_{k = 0}^{n} \frac{\partial^n M}{{\partial x}^k {\partial y}^{n - k}}(x , x) \, T^k .
$$
Then the $F_n$'s are reciprocal polynomials in $T$ and satisfy the recursive relation:
$$
F_n(x , T) = \left(\frac{T}{T + 1}\right) \dfrac{\partial F_{n - 1}}{\partial x}(x , T) + \dfrac{\partial^n M}{{\partial x}^n}(x , x) \cdot \frac{T^{n + 1} + 1}{T + 1} \quad (\forall n \in \N) .
$$
\end{prop} 

\begin{proof}
The fact that $F_n$ ($n \in \N_0$) is a reciprocal polynomial in $T$ follows from the symmetry of $M$. Next, for all $n \in \N$, we have
\begin{align*}
\dfrac{\partial F_{n - 1}}{\partial x}(x ,T) & = \sum_{k = 0}^{n - 1} \dfrac{d}{d x}\left(\dfrac{\partial^{n - 1} M}{{\partial x}^k {\partial y}^{n - 1 - k}}(x , x)\right) T^k \\
& \hspace*{-1cm} = \sum_{k = 0}^{n - 1} \left(\dfrac{\partial^n M}{{\partial x}^{k + 1} {\partial y}^{n - 1 - k}}(x , x) + \dfrac{\partial^n M}{{\partial x}^k {\partial y}^{n - k}}(x , x)\right) T^k \\
& \hspace*{-1cm} = \sum_{k = 0}^{n - 1} \dfrac{\partial^n M}{{\partial x}^{k + 1} {\partial y}^{n - 1 - k}}(x , x) \, T^k + \sum_{k = 0}^{n - 1} \dfrac{\partial^n M}{{\partial x}^k {\partial y}^{n - k}}(x , x) \, T^k \\
& \hspace*{-1cm} = \sum_{k = 1}^{n} \dfrac{\partial^n M}{{\partial x}^k {\partial y}^{n - k}}(x , x) \, T^{k - 1} + \sum_{k = 0}^{n - 1} \dfrac{\partial^n M}{{\partial x}^k {\partial y}^{n - k}}(x , x) \, T^k \\
& \hspace*{-1cm} = \sum_{k = 0}^{n} \dfrac{\partial^n M}{{\partial x}^k {\partial y}^{n - k}}(x , x) \, T^{k - 1} + \sum_{k = 0}^{n} \dfrac{\partial^n M}{{\partial x}^k {\partial y}^{n - k}}(x , x) \, T^k - \dfrac{\partial^n M}{{\partial y}^n}(x , x) T^{-1} - \dfrac{\partial^n M}{{\partial x}^n}(x , x) T^n \\
& \hspace*{-1cm} = \left(\frac{1}{T} + 1\right) F_n(x , T) - \dfrac{\partial^n M}{{\partial x}^n}(x , x) \left(T^n + \frac{1}{T}\right)
\end{align*}
(since $\frac{\partial^n M}{{\partial x}^n}(x , x) = \frac{\partial^n M}{{\partial y}^n}(x , x)$, due to the symmetry of $M$). Rearranging and simplifying yields the required recursive relation, completing the proof.
\end{proof}

Given a mean $M$ that admits an expansion near the first bissector, the following proposition provides an explicit formula expressing the functions $f_n$ in Formula \eqref{eq2} as linear combinations of the partial derivatives of $M$ of order $2 n$ at points on the first bisector, with rational coefficients that are independent of $M$.

\begin{prop}\label{p4}
Let $M$ be a mean, and for all $n \in \N$, let $f_n$ be the function given by Formula \eqref{eq1}. Then, for all $n \in \N$, we have:
\begin{equation}\label{eq9}
f_n(x) = \frac{1}{(2 n)! 2^{2 n - 1}} \sum_{k = 0}^{n - 1} (-1)^k \binom{2 n}{k} \dfrac{\partial^{2 n} M}{{\partial x}^k {\partial y}^{2 n - k}}(x , x) + \frac{(-1)^n}{n!^2 4^n} \dfrac{\partial^{2 n} M}{{\partial x}^n {\partial y}^n}(x , x) . 
\end{equation}
\end{prop}

\begin{proof}
For $x \in (0 , + \infty)$, let $\varphi_x$ be the function defined in \eqref{eq3}, namely $\varphi_x(t) := M(x + \frac{t}{2} , x - \frac{t}{2})$ for all $t \in (- 2 x , 2 x)$. A straightforward induction shows that for all $n \in \N_0$:
$$
\varphi_x^{(n)}(t) = \frac{1}{2^n} \sum_{k = 0}^{n} (-1)^{n - k} \binom{n}{k} \dfrac{\partial^n M}{{\partial x}^k {\partial y}^{n - k}}\left(x + \frac{t}{2} , x - \frac{t}{2}\right) .
$$
Using this formula in \eqref{eq1}, we deduce that for all $n \in \N$:
$$
f_n(x) := \dfrac{\varphi_x^{(2 n)}(0)}{(2 n)!} = \frac{1}{(2 n)! 4^n} \sum_{k = 0}^{2 n} (-1)^k \binom{2 n}{k} \dfrac{\partial^{2 n} M}{{\partial x}^k {\partial y}^{2 n - k}}(x , x) .
$$
Grouping together each term of order $k$ ($0 \leq k \leq n - 1$) with its symmetric counterpart of order $(2 n - k)$, which are equal due to the symmetry of $M$, we obtain:
\begin{align*}
f_n(x) & = \frac{2}{(2 n)! 4^n} \sum_{k = 0}^{n - 1} (-1)^k \binom{2 n}{k} \dfrac{\partial^{2 n} M}{{\partial x}^k {\partial y}^{2 n - k}}(x , x) + \frac{1}{(2 n)! 4^n} (-1)^n \binom{2 n}{n} \dfrac{\partial^{2 n} M}{{\partial x}^n {\partial y}^n}(x , x) \\
& = \frac{1}{(2 n)! 2^{2 n -1}} \sum_{k = 0}^{n - 1} (-1)^k \binom{2 n}{k} \dfrac{\partial^{2 n} M}{{\partial x}^k {\partial y}^{2 n - k}}(x , x) + \frac{(-1)^n}{n!^2 4^n} \dfrac{\partial^{2 n} M}{{\partial x}^n {\partial y}^n}(x , x) ,
\end{align*}
as required.
\end{proof}

From Propositions \ref{p3} and \ref{p4}, we immediately derive the following corollary:

\begin{coll}\label{c4}
Let $M$ be a mean, and for all $n \in \N$, let $f_n$ be the function given in Formula \eqref{eq1}. Then, for all $n \in \N$, $f_n$ can be expressed as a linear combination (with rational coefficients independent of $M$) of the characteristic functions $Q_M^{[i]}$ of $M$ ($1 \leq i \leq n$) and their higher-order derivatives. More precisely, for all $n \in \N$, there exist rational numbers $\alpha_i$ ($1 \leq i \leq n$), independent of $M$, such that the following identity holds:
\begin{equation}
f_n = \sum_{i = 1}^{n} \alpha_i \left(\frac{d}{d x}\right)^{\!\! 2 n - 2 i} \!\! Q_M^{[i]} . \tag*{$\square$}
\end{equation} 
\end{coll}

\begin{rmk}
Given a mean $M$, a closed-form expression for the rational coefficients in Corollary \ref{c4} appears difficult to determine. Instead, we propose to compute these coefficients for small values of $n$. Applying Formula \eqref{eq9} from Proposition \ref{p4} for $n = 1$ and using Formulas \eqref{eq-qay-2} and \eqref{eq-qay-3}, we obtain:
$$
f_1(x) = \frac{1}{4} \dfrac{\partial^2 M}{{\partial y}^2}(x , x) - \frac{1}{4} \dfrac{\partial^2 M}{\partial x \partial y}(x , x) = \frac{1}{2} Q_M^{[1]}(x) .
$$
Thus,
\begin{equation}\label{eq10}
f_1 = \frac{1}{2} Q_M^{[1]} .
\end{equation}
Applying Formula \eqref{eq9} for $n = 2$ and using Formulas \eqref{eq-qay-7} and \eqref{eq-qay-9}, we obtain
\begin{align*}
f_2(x) & = \frac{1}{192} \left(\dfrac{\partial^4 M}{{\partial y}^4}(x , x) - 4 \dfrac{\partial^4 M}{\partial x {\partial y}^3}(x , x)\right) + \frac{1}{64} \dfrac{\partial^4 M}{{\partial x}^2 {\partial y}^2}(x , x) \\
& = \left(\frac{1}{24} Q_M^{[2]} - \frac{1}{16} {Q_M^{[1]}}''\right)(x) .
\end{align*} 
Thus,
\begin{equation}\label{eq11}
f_2 = \frac{1}{24} Q_M^{[2]} - \frac{1}{16} {Q_M^{[1]}}'' .
\end{equation} 
\end{rmk}

Now, suppose that two given means $M_1$ and $M_2$ admit expansions near the first bisector, given by:
\begin{align*}
M_1(x , y) & = A + f_1(A) (x - y)^2 + f_2(A) (x - y)^4 + \dots , \\
M_2(x , y) & = A + g_1(A) (x - y)^2 + g_2(A) (x - y)^4 + \dots .
\end{align*}
Then, near the first bisector, we have:
\begin{align*}
M_1(x , y) - M_2(x , y) & = (f_1 - g_1)(A) (x - y)^2 + (f_2 - g_2)(A) (x - y)^4 + \dots \\
& = O\left((x - y)^2\right) . 
\end{align*}
For the difference between $M_1$ and $M_2$ to be of a smaller order than $O\left((x - y)^2\right)$, we must necessarily have $f_1 = g_1$, which, by virtue of \eqref{eq10}, is equivalent to stating that $Q_{M_1} = Q_{M_2}$. This motivates the following definitions:

\begin{defi}\label{defi2}
Let $M_1$ and $M_2$ be two means. We say that $M_1$ is \emph{close to $M_2$ near the first bisector} (or simply that $M_1$ is \emph{close to} $M_2$) if $M_1$ and $M_2$ have the same characteristic function; that is, if $Q_{M_1} = Q_{M_2}$.
\end{defi}

\begin{defi}\label{defi3}
A mean $M$ is said to be \emph{almost arithmetic} if it is close to the arithmetic mean; that is, if $Q_M = 0$ (since $Q_A = 0$). 
\end{defi}

\begin{expl}\label{expl1}
For all $x , y > 0$, define
$$
M(x , y) := \frac{x + y}{2} + \frac{1}{2} \cdot \dfrac{(x - y)^4}{1 + (x - y)^4} .
$$
It is clear that $M$ is infinitely differentiable and symmetric on $(0 , + \infty)^2$. Moreover, by distinguishing the two cases $\vert{x - y}\vert \leq 1$ and $\vert{x - y}\vert > 1$, we easily verify that for all $x , y > 0$,
$$
\dfrac{(x - y)^4}{1 + (x - y)^4} \leq \vert{x - y}\vert ,
$$
which implies
$$
\frac{x + y}{2} \leq M(x , y) \leq \max(x , y) .
$$
Thus, $M$ is a mean. Furthermore, for all sufficiently close $x , y > 0$, we have:
\begin{align*}
M(x , y) & = A + \frac{1}{2} (x - y)^4 \left(1 - (x - y)^4 + (x - y)^8 - \dots\right) \\[1mm]
& = A + \frac{1}{2} (x - y)^4 - \frac{1}{2} (x - y)^8 + \frac{1}{2} (x - y)^{12} - \dots ,
\end{align*}
showing that $M$ admits an expansion near the first bisector and is almost arithmetic.
\end{expl}

Note that this example is considered as a trivial case of an almost arithmetic mean. Nontrivial examples will be provided later (see \S \ref{subsubsec2-a}).

\subsubsection{The case of homogeneous and additively homogeneous means}

For means $M$ admitting an expansion near the first bisector, Corollaries \ref{c1} and \ref{c2} ensure that $Q_M(x)$ has the form $\frac{c}{x}$ for some $c \in \R$ if $M$ is homogeneous, and is constant if $M$ is additively homogeneous. The following proposition extends these facts to any mean.

\begin{prop}\label{p23}
Let $M$ be a mean. The following properties hold:
\begin{enumerate}
\item If $M$ is homogeneous, then its characteristic function has the form
$$
Q_M(x) = \dfrac{c}{x}
$$
for some constant $c \in \R$.
\item If $M$ is additively homogeneous, then its characteristic function is constant.
\end{enumerate}
\end{prop}

\begin{proof}~\\
{\large\bf 1.} Suppose that $M$ is homogeneous. Then, for all $\lambda , x , y > 0$, we have:
$$
M(\lambda x , \lambda y) = \lambda M(x , y) .
$$
Differentiating twice with respect to $x$ and then setting $y = x$, we obtain:
$$
\lambda^2 Q_M(\lambda x) = \lambda Q_M(x) \quad (\forall \lambda , x > 0) .
$$
That is,
$$
\lambda x Q_M(\lambda x) = x Q_M(x) \quad (\forall \lambda , x > 0) ,
$$
which implies that the function $x \mapsto x Q_M(x)$ is constant. Therefore, $Q_M(x) = \frac{c}{x}$ for some constant $c \in \R$, as required. \\[1mm]
{\large\bf 2.} Now, suppose that $M$ is additively homogeneous. Then, for all $\lambda , x, y > 0$, we have:
$$
M(x + \lambda , y + \lambda) = M(x , y) + \lambda .
$$
Differentiating twice with respect to $x$ and then setting $y = x$, we get
$$
Q_M(x + \lambda) = Q_M(x) \quad (\forall \lambda , x > 0) ,
$$
which shows that $Q_M$ is constant, as required. This completes the proof.
\end{proof}

\subsubsection{The characteristic functions of some common means and applications}\label{subsubsec2-a}

By direct calculation or by using Corollary \ref{c3} together with Formula \eqref{eq10}, we establish that the characteristic functions of the arithmetic mean $A$, the geometric mean $G$, the harmonic mean $H$, and the exponential mean $E$ are given by:
\begin{align*}
Q_A(x) & = 0 , \\
Q_G(x) & = - \frac{1}{4 x} , \\[1mm]
Q_H(x) & = - \frac{1}{2 x} , \\[1mm]
Q_E(x) & = \frac{1}{12} .
\end{align*}
Further calculations yield:
\begin{align*}
Q_{\eta_{a , b}}(x) & = \frac{a + b - 1}{4} \cdot \frac{1}{x} \quad (\forall a , b \in \R) , \\[2mm]
Q_{S_p}(x) & = \frac{p - 2}{12} \cdot \frac{1}{x} \quad (\forall p \in \R) .
\end{align*}
Using Formula \eqref{eq-intro-1}, we derive the following characteristic functions for the power means and the Lehmer means:
\begin{align}
Q_{\pi_r}(x) & = \frac{r - 1}{4} \cdot \frac{1}{x} \quad (\forall r \in \R) , \label{eq12} \\[1mm]
Q_{\L_t}(x) & = \frac{t - 1}{2} \cdot \frac{1}{x} \quad (\forall t \in \R) . \label{eq13}
\end{align} 
Next, substituting $p = 0$ and $p = 1$ into the above formula for $Q_{S_p}(x)$ yields
\begin{align*}
Q_L(x) & = - \frac{1}{6 x} , \\[2mm]
Q_I(x) & = - \frac{1}{12 x} .
\end{align*}
These examples show that the Gini mean $\eta_{a , b}$ ($a , b \in \R$) is almost arithmetic if and only if $a + b = 1$. This leads to the following notable approximation near the first bisector:
$$
\left(\dfrac{x^a + y^a}{x^{1 - a} + y^{1 - a}}\right)^{\frac{1}{2 a - 1}} \simeq \dfrac{x + y}{2} ,
$$
which is valid for all $a \in \R \setminus \{\frac{1}{2}\}$ with an error of order $O((x - y)^4)$). According to \eqref{eq-intro-2}, when $a = \frac{1}{2}$, this approximation must be replaced by:
$$
\exp\left(\dfrac{(\log x) \sqrt{x} + (\log y) \sqrt{y}}{\sqrt{x} + \sqrt{y}}\right) \simeq \dfrac{x + y}{2} .
$$
Furthermore, the Gini mean $\eta_{a , b}$ is close to the geometric mean $G$ if and only if $b = - a$, which means that it coincides with $G$ (as it is straightforward to verify that $\eta_{a , - a} = G$ for all $a \in \R$). Moreover, two power means are close if and only if they are identical, and the same holds for two Lehmer means.

Another application consists in observing that a power mean $\pi_r$ is close to a Lehmer mean $\L_t$ if and only if $r = 2 t - 1$. This leads to the following notable approximation near the first bisector:
$$
\left(\dfrac{x^{2 t - 1} + y^{2 t - 1}}{2}\right)^{\frac{1}{2 t -1}} \simeq \dfrac{x^t + y^t}{x^{t - 1} + y^{t - 1}} ,
$$
which is valid for all $t \in \R \setminus \{\frac{1}{2}\}$ with an error of order $O((x - y)^4)$. Note that this approximation becomes an equality for $t \in \{0 , 1\}$, as well as in the limit $t \to \frac{1}{2}$.

We now determine the characteristic function of the Gauss Arithmetic-Geometric Mean ($\AGM$). To do so, we rely on the following proposition:

\begin{prop}\label{p5}
Let $M_0 , M_1 , M_2$ be three means, and let $M$ be the mean defined by:
$$
M := M_0\left(M_1 , M_2\right) .
$$
Then, the characteristic function of $M$ is given by:
$$
Q_M = \dfrac{Q_{M_1} + Q_{M_2}}{2} ,
$$
which is independent of $M_0$.
\end{prop}

\begin{proof}
A direct computation involving the partial derivatives of composite and product functions in two variables yields:
\begin{multline*}
\dfrac{\partial^2 M}{{\partial x}^2} = \dfrac{\partial^2 M_1}{{\partial x}^2} \cdot \dfrac{\partial M_0}{\partial x}\left(M_1 , M_2\right) + \dfrac{\partial^2 M_2}{{\partial x}^2} \cdot \dfrac{\partial M_0}{\partial y}\left(M_1 , M_2\right) + \left(\dfrac{\partial M_1}{\partial x}\right)^2 \cdot \dfrac{\partial^2 M_0}{{\partial x}^2}\left(M_1 , M_2\right) \\
+ \left(\dfrac{\partial M_2}{\partial x}\right)^2 \cdot \dfrac{\partial^2 M_0}{{\partial y}^2}\left(M_1 , M_2\right) + 2 \, \dfrac{\partial M_1}{\partial x} \cdot \dfrac{\partial M_2}{\partial x} \cdot \dfrac{\partial^2 M_0}{\partial x \partial y}\left(M_1 , M_2\right) .
\end{multline*}
Evaluating at $(x , x)$ ($x > 0$), we obtain:
\begin{multline*}
Q_M(x) := \dfrac{\partial^2 M}{{\partial x}^2}(x , x) = \dfrac{\partial^2 M_1}{{\partial x}^2}(x , x) \cdot \dfrac{\partial M_0}{\partial x}(x , x) + \dfrac{\partial^2 M_2}{{\partial x}^2}(x , x) \cdot \dfrac{\partial M_0}{\partial y}(x , x) \\
+ \left(\dfrac{\partial M_1}{\partial x}(x , x)\right)^2 \cdot \dfrac{\partial^2 M_0}{{\partial x}^2}(x , x) + \left(\dfrac{\partial M_2}{\partial x}(x , x)\right)^2 \cdot \dfrac{\partial^2 M_0}{{\partial y}^2}(x , x) \\
+ 2 \, \dfrac{\partial M_1}{\partial x}(x , x) \cdot \dfrac{\partial M_2}{\partial x}(x , x) \cdot \dfrac{\partial^2 M_0}{\partial x \partial y}(x , x) .
\end{multline*}
Using Formulas \eqref{eq-qay-1}, \eqref{eq-qay-2}, and \eqref{eq-qay-3}, the result follows.
\end{proof}

\begin{rmk}
If the means $M_0 , M_1$, and $M_2$ in Proposition \ref{p5} admit expansions near the first bisector, an alternative proof (involving simpler calculations) can be derived using these expansions. 
\end{rmk}

Applying Proposition \ref{p5}, we obtain the following corollary:

\begin{coll}\label{c5}
The characteristic function of the Gauss Arithmetic-Geometric Mean is given by:
$$
Q_{\AGM}(x) = - \frac{1}{8 x} ,
$$
for all $x > 0$.
\end{coll}

\begin{proof}
By definition, the Arithmetic-Geometric Mean satisfies
$$
\AGM = \AGM\left(A , G\right) ,
$$
and applying Proposition \ref{p5}) gives
$$
Q_{\AGM} = \dfrac{Q_A + Q_G}{2} .
$$
Since we previously established that $Q_A(x) = 0$ and $Q_G(x) = - \frac{1}{4 x}$ for all $x > 0$, the result follows immediately.
\end{proof}

\begin{rmk}
The characteristic function of the Gauss Arithmetic-Geometric Mean can also be computed from Formula \eqref{eq7}.
\end{rmk}

\subsubsection{Comparison of two means using their characteristic functions}\label{subsubsec2-b}

As seen above, the characteristic function of a mean measures its proximity to the arithmetic mean near the first bisector. Therefore, comparing two means (in a neighborhood of the first bisector) reduces to comparing their characteristic functions, which are functions of a single variable. This significantly simplifies the process of deriving various inequalities between the means, as demonstrated below. We state the following proposition:

\begin{prop}\label{p22}
let $M_1$ and $M_2$ be two means satisfying:
$$
Q_{M_1}(x) > Q_{M_2}(x) \quad (\forall x \in (0 , + \infty)) .
$$
Then we have, near the first bisector:
$$
M_1(x , y) \geq M_2(x , y) .
$$
\end{prop}

\begin{proof}
For each $x > 0$, define the function $\Phi_x : (- 2 x , 2 x) \rightarrow \R$ by
$$
\Phi_x(t) := M_1\left(x + \frac{t}{2} , x - \frac{t}{2}\right) - M_2\left(x + \frac{t}{2} , x - \frac{t}{2}\right) \quad (\forall t \in (- 2 x , 2 x)) .
$$
A straightforward computation gives for all $x > 0$:
$$
\Phi_x''(0) = Q_{M_1}(x) - Q_{M_2}(x) > 0 .
$$
Hence, for all $x > 0$, the function $\Phi_x$ is convex near $t = 0$. It follows that for all $x > 0$ and all $t$ sufficiently close to $0$, we have
$$
\Phi_x(t) \geq \Phi_x(0) + \Phi_x'(0) t ;
$$
that is,
\begin{equation}\label{eq26}
\Phi_x(t) \geq 0 
\end{equation}
(since $\Phi_x(0) = \Phi_x'(0) = 0$). Finally, for $x , y > 0$, sufficiently close, substituting in \eqref{eq26} $x$ by $\frac{x + y}{2}$ and $t$ by $(x - y)$ yields
$$
M_1(x , y) - M_2(x , y) \geq 0 ,
$$
as required.
\end{proof}

\begin{coll}\label{c16}
Let $M_1$ and $M_2$ be two homogeneous means satisfying:
$$
Q_{M_1}(x) > Q_{M_2}(x) \quad (\forall x \in (0 , + \infty)) .
$$
Then, for all $x , y > 0$, we have:
$$
M_1(x , y) \geq M_2(x , y) .
$$
\end{coll}

\begin{proof}
By Proposition \ref{p22}, the inequality $M_1(x , y) \geq M_2(x , y)$ holds in a neighborhood of the first bisector. The homogeneity of $M_1$ and $M_2$ then implies that the inequality holds everywhere on $(0 , + \infty)^2$.
\end{proof}

\begin{expls}
From \S \ref{subsubsec2-a}, we have, for all $x > 0$:
$$
Q_A(x) > Q_I(x) > Q_{\AGM}(x) > Q_L(x) > Q_G(x) > Q_H(x) .
$$
By Corollary \ref{c16}, we obtain the chain of inequalities between the corresponding means:
$$
A \geq I \geq \AGM \geq L \geq G \geq H .
$$
Regarding parameter means, we also refer to \S \ref{subsubsec2-a}, from which it follows that for all $x > 0$:
\begin{align*}
Q_{\pi_r}(x) & > Q_{\pi_s}(x) \quad\quad (\text{whenever } r > s) , \\
Q_{\L_t}(x) & > Q_{\L_s}(x) \quad\quad (\text{whenever } t > s) , \\
Q_{S_p}(x) & > Q_{S_q}(x) \quad\quad (\text{whenever } p > q) , \\
Q_{\eta_{a , b}}(x) & > Q_{\eta_{c , d}}(x) \quad\quad (\text{whenever } a + b > c + d) .
\end{align*}
Applying again Corollary \ref{c16}, we deduce the corresponding inequalities between means:
\begin{align*}
\pi_r & \geq \pi_s \quad\quad (\text{whenever } r > s) , \\
\L_t & \geq \L_s \quad\quad (\text{whenever } t > s) , \\
S_p & \geq S_q \quad\quad (\text{whenever } p > q) , \\
\eta_{a , b} & \geq \eta_{c , d} \quad\quad (\text{whenever } a + b > c + d) .
\end{align*}
Note that all of these inequalities are already known in the literature. For instance, the last one is contained in a result by Elezovi\`c and Vuk\v{s}i\'{c} \cite[Theorem A]{elv2}.
\end{expls}

\section{Classification of means}

In what follows, we study various classes of means; some of these have already been previously considered by earlier authors, while others are new. For each of these classes, we determine the form of the corresponding characteristic function, which sometimes allows us to identify intersections between different classes of means. Furthermore, each of the classes below satisfies two important properties, which we refer to as \emph{minimality} and \emph{completeness}.

\begin{defis}\label{defis1}
Let $\mathscr{C}$ be a class of means. 
\begin{enumerate}
\item We say that $\mathscr{C}$ is \emph{minimal} if for all $M_1 , M_2 \in \mathscr{C}$, we have:
$$
Q_{M_1} = Q_{M_2} \Longrightarrow M_1 = M_2 .
$$
In other words, $\mathscr{C}$ is minimal if any two means in $\mathscr{C}$ that are close must be identical.
\item We say that $\mathscr{C}$ is \emph{complete} if for every real-valued infinitely differentiable function $f$ on $(0 , + \infty)$, there exists a mean $M$ in $\mathscr{C}$ such that $Q_M = f$. 
\end{enumerate}
\end{defis}

\begin{rmk}
Let $\mathscr{C}$ be a class of means and define the mapping
$$
\begin{array}{rcl}
\Phi :~ \mathscr{C} & \longrightarrow & C^{\infty}\left((0 , + \infty) , \R\right) \\
M & \longmapsto & Q_M
\end{array} .
$$
Then, saying that $\mathscr{C}$ is minimal is equivalent to $\Phi$ being injective, and saying that $\mathscr{C}$ is complete is equivalent to $\Phi$ being surjective. Consequently, saying that $\mathscr{C}$ is both minimal and complete is equivalent to $\Phi$ being bijective. Thus, in a minimal and complete class of means, each mean is uniquely determined by its characteristic function. 
\end{rmk}

\subsection{Normal means}

This class was introduced for the first time by the author \cite{far}, although the definition of a mean therein differs slightly from that of the present paper.

\begin{defi}[\cite{far}]\label{defi4}
A mean $M$ is said to be \emph{normal} if it has the form:
$$
M(x , y) = \dfrac{x P(x) + y P(y)}{P(x) + P(y)} \quad (\forall x , y > 0) ,
$$
where $P$ is a positive infinitely differentiable function on $(0 , + \infty)$, called the \emph{weight function} associated with $M$, which is uniquely defined up to a multiplicative positive constant.
\end{defi}

The class of normal means includes the Lehmer means. More precisely, every Lehmer mean $\L_t$ ($t \in \R$) is a normal mean with the associated weight function $P(x) = x^{t - 1}$. In particular, the arithmetic mean $A$, the geometric mean $G$, and the harmonic mean $H$ are all normal means.

The following proposition provides the expression of the characteristic function of a normal mean in terms of its associated weight function.

\begin{prop}\label{p6}
Let $P$ be a positive infinitely differentiable function on $(0 , + \infty)$, and let $M$ be the associated normal mean. Then, we have
$$
Q_M = \frac{1}{2} \frac{P'}{P} .
$$
\end{prop} 

\begin{proof}
To simplify the calculation, observe that for all $x , y > 0$, we have:
$$
M(x , y) = x + (y - x) \dfrac{P(y)}{P(x) + P(y)} .
$$
Using this expression, we compute:
$$
\dfrac{\partial^2 M}{{\partial x}^2}(x , y) = \dfrac{2 P'(x) P(y)}{\left(P(x) + P(y)\right)^2} + (x - y) P(y) \left[\dfrac{P''(x) \left(P(x) + P(y)\right) - 2 P'^2(x)}{\left(P(x) + P(y)\right)^3}\right] .
$$
The result follows by setting $y = x$.
\end{proof}

From Proposition \ref{p6}, we derive the following corollary:

\begin{coll}\label{c6}
The class of normal means is minimal and complete.
\end{coll}

\begin{proof}
We first show that the class of normal means is minimal. Let $M_1$ and $M_2$ be two normal means with the same characteristic function, and let $P_1$ and $P_2$ be their respective associated weight functions. By Proposition \ref{p6}, we have
$$
\frac{1}{2} \dfrac{P_1'}{P_1} = \frac{1}{2} \dfrac{P_2'}{P_2} ,
$$
implying that there exists a real constant $c$ (necessarily positive) such that $P_2 = c P_1$. Thus, $M_1 = M_2$, as required.

Next, we prove completeness. Let $f \in C^{\infty}\left((0 , + \infty) , \R\right)$ be arbitrary. Defining $F$ as an antiderivative of $f$ and setting $P := e^{2 F}$, which is a positive infinitely differentiable function on $(0 , + \infty)$, the normal mean $M$ having $P$ as a weight function satisfies (according to Proposition \ref{p6}):
$$
Q_M = \frac{1}{2} \dfrac{P'}{P} = f .
$$
This confirms that normal means form a complete class of means. The proof is complete.
\end{proof}

As an application of the characteristic functions of normal means, we establish the following result:

\begin{thm}\label{t2}
The only power means that are normal are the arithmetic mean, the geometric mean, and the harmonic mean.
\end{thm}

\begin{proof}
We already know that the means $A , G$, and $H$ are both power means and normal means. Conversely, let $M$ be a normal power mean and assume, for contradiction, that $M$ is neither $A$, $G$, nor $H$. Let $P$ be the weight function associated with $M$ (as a normal mean), and let $r \in \R \setminus \{-1 , 0 , 1\}$ such that $M = \pi_r$. Then, on the one hand, we have $Q_M = \frac{1}{2} \frac{P'}{P}$ by proposition \ref{p6}, and on the other hand, from \eqref{eq12}, we have $Q_M(x) = \frac{r - 1}{4} \frac{1}{x}$. Equating these two expressions of $Q_M$ gives the differential equation
$$
\frac{1}{2} \dfrac{P'}{P} = \frac{r - 1}{4} \frac{1}{x} .
$$
Solving this, we obtain
$$
P(x) = c \, x^{\frac{r - 1}{2}} ,
$$
where $c$ is a positive constant. Thus, we have
$$
M(x , y) = \dfrac{x P(x) + y P(y)}{P(x) + P(y)} = \dfrac{x^{\frac{r + 1}{2}} + y^{\frac{r + 1}{2}}}{x^{\frac{r - 1}{2}} + y^{\frac{r - 1}{2}}} \quad (\forall x , y > 0) ,
$$
and also
$$
M(x , y) = \pi_r(x , y) = \left(\dfrac{x^r + y^r}{2}\right)^{1/r} \quad (\forall x , y > 0) .
$$
Equating these two expressions of $M(x , y)$, rearranging, and setting $y = 1$ leads to
\begin{equation}\label{eq14}
\left(\dfrac{x^{\frac{r + 1}{2}} + 1}{x^{\frac{r - 1}{2}} + 1}\right)^r = \dfrac{x^r + 1}{2} \quad (\forall x > 0) . 
\end{equation}
To show that \eqref{eq14} is impossible, we consider the following cases: \\[1mm]
\textbf{Case 1:} (if $r > 1$). In this case, the left-hand side of \eqref{eq14} is asymptotically $x^r$ as $x \to + \infty$, while the right-hand side is asymptotically $\frac{1}{2} x^r$. Thus, \eqref{eq14} is impossible for this case. \\[1mm]
\textbf{Case 2:} (if $0 < r < 1$). In this case, the left-hand side of \eqref{eq14} is asymptotically $x^{\frac{r (r + 1)}{2}}$, while the right-hand side is asymptotically $\frac{1}{2} x^r$ as $x \to + \infty$. Since $\frac{r (r + 1)}{2} < r$ then \eqref{eq14} is also impossible for this case. \\[1mm]
\textbf{Case 3:} (if $-1 < r < 0$). In this case, the left-hand side of \eqref{eq14} tends to $0$, while the right-hand side tends to $\frac{1}{2}$ as $x \to + \infty$. Thus \eqref{eq14} is also impossible for this case. \\[1mm]
\textbf{Case 4:} (if $r < -1$). In this case, the left-hand side of \eqref{eq14} tends to $1$, while the right-hand side tends to $\frac{1}{2}$ as $x \to + \infty$. Thus \eqref{eq14} is also impossible for this case. 

In conclusion, \eqref{eq14} is impossible, confirming that $M$ must be either the arithmetic, geometric, or harmonic mean. This completes the proof. 
\end{proof}

\begin{rmk}
Using Proposition \ref{p23}, we fairly easily show (as done in the proof of Theorem \ref{t2}) that the only homogeneous means that are normal are the Lehmer means.
\end{rmk}

\subsubsection*{Characterization of strict normal means}
\addcontentsline{toc}{subsubsection}{Characterization of strict normal means}

If $M$ is a strict mean (see \S\ref{sec1}), we can associate to $M$ (as was done in \cite{far}) the important function $\varphi_M : (0 , + \infty)^2 \rightarrow (0 , + \infty)$, defined by:
\begin{equation}\label{eq15}
\varphi_M(x , y) := \begin{cases}
\log\left(- \dfrac{M(x , y) - x}{M(x , y) - y}\right) & \text{if } x \neq y \\
0 & \text{if } x = y
\end{cases} \quad (\forall x , y > 0) .
\end{equation}
We can verify that $\varphi_M$ is well-defined, infinitely differentiable on $(0 , + \infty)^2$, and asymmetric; namely, $\varphi_M(x , y) = - \varphi_M(y , x)$ for all $x , y > 0$ (see \cite[Theorem 2.1]{far}). In particular, we have $\varphi_M(x , x) = 0$ for all $x > 0$. Interestingly, if $M$ is a strict normal mean associated with some weight function $P$, then the expression of $\varphi_M$ simplifies to
$$
\varphi_M(x , y) = \log{P(x)} - \log{P(y)}  \quad (\forall x , y > 0) ,
$$
which is a sum of two functions, one depending only on $x$ and the other only on $y$. Conversely, if a strict mean $M$ has an associated function $\varphi_M$ of the form
\begin{equation}\label{eq16}
\varphi_M(x , y) = f(x) + g(y) ,
\end{equation}
for some functions $f , g$ on $(0 , + \infty)$, then, because $\varphi_M(x , x) = 0$, we must necessarily have $g = - f$. Hence, from \eqref{eq15} (see also \cite[Theorem 2.1]{far}), we derive
$$
M(x , y) = \dfrac{x e^{f(x)} + y e^{f(y)}}{e^{f(x)} + e^{f(y)}} \quad (\forall x , y > 0) ,
$$
showing that $M$ is the normal mean associated with the weight function $e^f$. It follows from this reasoning that the form \eqref{eq16} of $\varphi_M$ characterizes the normal strict means. Finally, since this particular form is clearly characterized by the partial differential equation $\frac{\partial^2 \varphi_M}{\partial x \partial y} = 0$, we conclude that a strict mean $M$ is normal if and only if $\frac{\partial^2 \varphi_M}{\partial x \partial y} = 0$. Developing this, we obtain the following proposition:

\begin{prop}\label{p7}
A strict mean $M$ is normal if and only if it satisfies the partial differential equation:
\begin{equation}
(x - y) \left[\dfrac{\partial^2 M}{\partial x \partial y} \cdot (M - x) \cdot (M - y) - 2 \, \dfrac{\partial M}{\partial x} \cdot \dfrac{\partial M}{\partial y} \cdot (M - A)\right] = \dfrac{\partial M}{\partial x} \cdot (M - x)^2 - \dfrac{\partial M}{\partial y} \cdot (M - y)^2 . \tag*{$\square$}
\end{equation}
\end{prop}

\begin{rmk}
Curiously, for the strict normal means $A , G$, and $H$, both sides of the partial differential equation in Proposition \ref{p7} are zero!
\end{rmk}

\begin{rmk}
In principle, Proposition \ref{p7} should provide us with an alternative proof of Theorem \ref{t2}. However, the calculations it involves are quite tedious.
\end{rmk}

\subsection{Additive means}\label{subsec3-2}

This class is well known in the mathematical literature under different names. In 
\cite[\S 4]{bul} for example, it is referred to as \emph{quasi-arithmetic means}.

\begin{defi}\label{defi5}
A mean $M$ is said to be \emph{additive} if it has the form:
$$
M(x , y) = f^{-1}\left(\dfrac{f(x) + f(y)}{2}\right) \quad (\forall x , y > 0) ,
$$
where $f$ is a strictly monotonic and infinitely differentiable function on $(0 , + \infty)$. Such a function $f$ is called \emph{an additivity function} associated with $M$. 
\end{defi}

The class of additive means includes power means. More precisely, every power mean $\pi_r$ ($r \in \R$) is an additive mean with an associated additivity function $f$ given by
$$
f(x) = \begin{cases}
x^r & \text{if } r \neq 0 \\
\log{x} & \text{if } r = 0
\end{cases} \quad (\forall x > 0) .
$$
In particular, the arithmetic mean $A$, the geometric mean $G$, and the harmonic mean $H$ are all additive means.

It is straightforward to verify that if $f$ is an additivity function associated with an additive mean $M$, and if $\alpha , \beta$ are two real numbers with $\alpha \neq 0$, then the function $(\alpha f + \beta)$ is also an additivity function associated with $M$. Note that the converse of this property is also true and can be deduced from a more general result given in \cite[\S 4]{bul}. More precisely, we have the following proposition:

\begin{prop}[Deduced from {\cite[\S 4 , Theorem 5]{bul}}]\label{p8}
If $f$ and $g$ are two additivity functions associated with a same additive mean $M$, then there exist $\alpha , \beta \in \R$, with $\alpha \neq 0$, such that $g = \alpha f + \beta$. \hfill $\square$
\end{prop}

\begin{rmk}
Using the concept of ``characteristic function of a mean'', we will be able to provide a short proof of Proposition \ref{p8} (see Remark \ref{r2} below).
\end{rmk}

We first establish the expression of the characteristic function of an additive mean in terms of an associated additivity function.

\begin{prop}\label{p9}
Let $f$ be a strictly monotonic and infinitely differentiable function on $(0 , + \infty)$, and let $M$ be the associated additive mean. Then, we have
$$
Q_M = \frac{1}{4} \dfrac{f''}{f'} .
$$
\end{prop}

\begin{proof}
From the expression of $M$ in terms of $f$, we have for all $x , y > 0$:
$$
f\left(M(x , y)\right) = \frac{1}{2} \left(f(x) + f(y)\right) .
$$
Differentiating twice with respect to $x$ yields
$$
\dfrac{\partial^2 M}{{\partial x}^2}(x , y) \cdot f'\left(M(x , y)\right) + \left(\dfrac{\partial M}{\partial x}(x , y)\right)^2 f''\left(M(x , y)\right) = \frac{1}{2} f''(x) .
$$
The result follows by setting $y = x$ and using Formula \eqref{eq-qay-1}.
\end{proof}

\begin{rmk}\label{r2}
Proposition \ref{p8} follows immediately from Proposition \ref{p9}. Indeed, if $f$ and $g$ are two additivity functions associated with a same additive mean $M$, then, according to Proposition \ref{p9}, we have 
$$
Q_M = \frac{1}{4} \frac{f''}{f'} = \frac{1}{4} \frac{g''}{g'} ,
$$
implying that the function $\frac{g'}{f'}$ is constant. Thus, there exist $\alpha , \beta \in \R$, with $\alpha \neq 0$, such that $g = \alpha f + \beta$, as required.
\end{rmk}

From Proposition \ref{p9}, we derive the following corollary:

\begin{coll}\label{c7}
The class of additive means is minimal and complete.
\end{coll}

\begin{proof}
We first show that the class of additive means is minimal. Let $M_1$ and $M_2$ be two additive means with the same characteristic function, and let $f_1$ and $f_2$ be their respective associated additivity functions. By proposition \ref{p9}, we have
$$
\frac{1}{4} \dfrac{f_1''}{f_1'} = \frac{1}{4} \dfrac{f_2''}{f_2'} ,
$$
implying that $f_1$ and $f_2$ are related by a relation of the form $f_2 = \alpha f_1 + \beta$ ($\alpha , \beta \in \R$, $\alpha \neq 0$). Thus, $M_1 = M_2$, as required.

Next, we prove completeness. Let $u \in C^{\infty}((0 , + \infty) , \R)$ be arbitrary. Defining $v$ as an antiderivative of $u$ and $w$ as an antiderivative of $e^{4 v}$, it is clear that $w$ is strictly increasing and infinitely differentiable on $(0 , + \infty)$ and that the additive mean $M$ having $w$ as an associated additivity function satisfies (according to Proposition \ref{p9}):
$$
Q_M = \frac{1}{4} \dfrac{w''}{w'} = u .
$$
This confirms that additive means form a complete class of means, completing the proof.
\end{proof}

As an application of the characteristic functions of additive means, we establish the following results:

\begin{thm}\label{t3}
The only additive means that are homogeneous are the power means.
\end{thm}

\begin{proof}
It is evident that power means are both additive and homogeneous. Conversely, let $M$ be a homogeneous additive mean. On the one hand, from Proposition \ref{p23}, we know that the characteristic function of $M$ has the form $Q_M(x) = \frac{c}{x}$ ($c \in \R$). On the other hand, denoting by $f$ an additivity function associated with $M$ (as an additive mean), we have (by Proposition \ref{p9}): $Q_M = \frac{1}{4} \frac{f''}{f'}$. Equating these two expressions for $Q_M$ yields the differential equation
$$
\frac{1}{4} \dfrac{f''}{f'} = \dfrac{c}{x} ,
$$
which integrates to
$$
f(x) = \begin{cases}
\alpha \, x^{4 c + 1} + \beta & \text{if } c \neq - \frac{1}{4} \\
\alpha \log{x} + \beta & \text{if } c = - \frac{1}{4} 
\end{cases}  \quad (\alpha , \beta \in \R , \alpha \neq 0) .
$$
Since $M$ remains unchanged when its additivity function $f$ is replaced by $\lambda f + \mu$ (for some $\lambda , \mu \in \R$, $\lambda \neq 0$), we can choose a simpler additivity function associated with $M$:
$$
\widetilde{f}(x) = \begin{cases}
x^{4 c + 1} & \text{if } c \neq - \frac{1}{4} \\
\log{x} & \text{if } c = - \frac{1}{4} 
\end{cases} .
$$
Consequently, $M$ is the power mean given by
$$
M = \pi_{4 c + 1} .
$$
This completes the proof.
\end{proof}

\begin{coll}\label{c8}
The only additive means that are Lehmer means are the arithmetic mean $A$, the geometric mean $G$, and the harmonic mean $H$. 
\end{coll}

\begin{proof}
It is well known that $A$, $G$, and $H$ are both additive and Lehmer means. Conversely, let $M$ be a Lehmer additive mean. Then $M$ is homogeneous and additive, which implies (by Theorem \ref{t3}) that $M$ is a power mean. Since $M$ is also normal (as it is Lehmer), we conclude (by Theorem \ref{t2}) that $M$ must be either $A$, $G$, or $H$, as required.
\end{proof}

We now turn our attention to means that are both normal and additive. For any mean $M$ and any nonnegative constant $c$, define the transformed mean $M_c$ by:
$$
M_c(x , y) := M(x + c , y + c) - c \quad (\forall x , y > 0) .
$$
We observe that for any $c \geq 0$, the classes of normal means and additive means remain invariant under the transformation $M \mapsto M_c$. More precisely, for any $c \geq 0$, the following properties hold:
\begin{itemize}
\item If $M$ is a normal mean with an associated weight function $P(x)$, then $M_c$ is a normal mean with the associated weight function $P(x + c)$.
\item If $M$ is an additive mean with an associated additivity function $f(x)$, then $M_c$ is an additive mean with the associated additivity function $f(x + c)$.
\end{itemize}
Since the means $A$, $G$, and $H$ are all normal and additive, it follows that for all $c \geq 0$, the means $A_c = A$, $G_c$, and $H_c$ remain all normal and additive. The explicit expressions for $G_c$ and $H_c$ are given by:
\begin{align}
G_c(x , y) & = \sqrt{(x + c) (y + c)} - c , \label{eq17} \\[1mm]
H_c(x , y) & = \dfrac{2 x y + c (x + y)}{x + y + 2 c} \label{eq18} 
\end{align}
(for all $x , y > 0$). To complement these simple findings, we pose the following open problem, which seems challenging to resolve using only the characteristic function of a mean.

\medskip

\noindent\textbf{Open problem.} \emph{Are the means $A$, $G_c$ ($c \geq 0$), and $H_c$ ($c \geq 0$) the only means that are both normal and additive?} 

\medskip

We end this subsection with another quite important property of additive means, stated in the following proposition:

\begin{prop}\label{p10}
Every additive mean is strictly isotone.
\end{prop}

\begin{proof}
Let $M$ be an additive mean with an associated additivity function $f$. Since $f$ and $f^{-1}$ 
share the same monotonicity (both being strictly increasing or both strictly decreasing), it follows that for all $y > 0$, the function
$$
x \longmapsto f^{-1}\left(\dfrac{f(x) + f(y)}{2}\right) = M(x , y)
$$
is strictly increasing on $(0 , + \infty)$; equivalently, $M$ is strictly isotone, as required.
\end{proof}

\begin{expl}
Using Proposition \ref{p10}, we can easily show that the Lehmer mean $\L_t$ ($t < 0$) is not additive, a fact already included in Corollary \ref{c8}. Let $t < 0$ and assume for contradiction that $\L_t$ is isotone. Then, for all $x \geq 2$:
$$
\L_t(x , 1) \geq \L_t(2 , 1) .
$$ 
Letting $x \to + \infty$ leads to
$$
1 \geq \L_t(2 , 1) ,
$$
which is a contradiction. Hence, $\L_t$ is not isotone, and by Proposition \ref{p10}, it is not additive.
\end{expl}

\subsection{Integral means of the first kind}\label{subsec3-3}

This class has been referred to by earlier authors as \emph{integral means} (see e.g., \cite[Chapter 8, Proposition 8.3]{bor}). However, to distinguish it from another similar class of means, we introduce the terms ``\emph{of the first kind}'' and ``\emph{of the second kind}'' to designate each of them respectively.

\begin{defi}\label{defi6}
A mean $M$ is called \emph{an integral mean of the first kind} if there exists a strictly monotonic and infinitely differentiable function $f$ on $(0 , + \infty)$ such that, for all $x , y > 0$ with $x \neq y$, we have:
\begin{equation}\label{eq19}
M(x , y) = f^{-1}\left(\frac{1}{y - x} \int_{x}^{y} f(t) \, d t\right) .
\end{equation}
Such a function $f$ is called \emph{an integrality function} associated with $M$.
\end{defi}

\begin{rmk}\label{r3}
By performing the change of variable $t = x u + y (1 - u)$ in the integral on the right-hand side of \eqref{eq19}, the formula becomes:
\begin{equation}\label{eq20}
M(x , y) = f^{-1}\left(\int_{0}^{1} f\left(x u + y (1 - u)\right) \, d u\right) ,
\end{equation}
which also holds when $y = x$, thereby showing that any continuous function on $(0 , + \infty)^2$ having the form in \eqref{eq19} (for $x \neq y$) defines a mean. More importantly, Formula \eqref{eq20} serves two further purposes:
\begin{enumerate}
\item[(1)] It allows an easy computation of the characteristic function of an integral mean of the first kind, as we will see below.
\item[(2)] It motivates the extension of the class of integral means of the first kind into a broader class of means (see \S \ref{subsec3-5}).
\end{enumerate} 
\end{rmk}

The class of integral means of the first kind includes the Stolarsky means. More precisely, every Stolarsky mean $S_p$ ($p \in \R$) is an integral mean of the first kind associated with the integrality function
$$
f(t) = \begin{cases}
t^{p - 1} & \text{if } p \neq 1 \\
\log{t} & \text{if } p = 1
\end{cases} \quad (\forall t > 0) .
$$
In particular, the arithmetic mean $A$, geometric mean $G$, logarithmic mean $L$, and identric mean $I$ all belong to this class. Additionally, the exponential mean $E$ is also an integral mean of the first kind, associated with the integrality function $f(t) = e^t$. However, it will be shown later that the harmonic mean $H$ does not belong to the class of integral means of the first kind.

It is easy to verify that if $f$ is an integrality function associated with an integral mean of the first kind $M$, then so is any affine transformation $\alpha f + \beta$ ($\alpha , \beta \in \R$, $\alpha \neq 0$). The converse will be established in Corollary \ref{c9} below.

We first express the characteristic function of an integral mean of the first kind in terms of its associated integrality function.

\begin{prop}\label{p11}
Let $f$ be a strictly monotonic and infinitely differentiable function on $(0 , + \infty)$, and let $M$ be the integral mean of the first kind associated with $f$. Then, we have
$$
Q_M = \frac{1}{12} \dfrac{f''}{f'} .
$$
\end{prop}

\begin{proof}
From Formula \eqref{eq20}, we have for all $x , y > 0$:
$$
f\left(M(x , y)\right) = \int_{0}^{1} f\left(x u + y (1 - u)\right) \, d u .
$$
Differentiating twice with respect to $x$ yields
$$
\dfrac{\partial^2 M}{{\partial x}^2}(x , y) f'\left(M(x , y)\right) + \left(\dfrac{\partial M}{\partial x}(x , y)\right)^2 f''\left(M(x , y)\right) = \int_{0}^{1} u^2 f''\left(x u + y (1 - u)\right) \, d u .
$$
The result follows by setting $y = x$ and using Formula \eqref{eq-qay-1}.
\end{proof}

From Proposition \ref{p11}, we derive the following corollaries:

\begin{coll}\label{c9}
If $f$ and $g$ are two integrality functions associated with a same integral mean of the first kind, then there exist $\alpha , \beta \in \R$, with $\alpha \neq 0$, such that $g = \alpha f + \beta$.
\end{coll}

\begin{proof}
Let $M$ be an integral mean of the first kind with two associated integrality functions $f$ and $g$. By proposition \ref{p11}, we have
$$
Q_M = \frac{1}{12} \dfrac{f''}{f'} = \frac{1}{12} \dfrac{g''}{g'} .
$$
It follows that $\frac{g'}{f'}$ is constant (nonzero), hence $g = \alpha f + \beta$ for some $(\alpha , \beta) \in \R^* \times \R$, as claimed.
\end{proof}

\begin{coll}\label{c10}
The class of integral means of the first kind is minimal and complete.
\end{coll}

\begin{proof}
\emph{Minimality}: Let $M_1$ and $M_2$ be two integral means of the first kind with the same characteristic function, and let $f_1$ and $f_2$ be their respective associated integrality functions. Then, by Proposition \ref{p11}, we have
$$
\frac{1}{12} \dfrac{f_1''}{f_1'} = \frac{1}{12} \dfrac{f_2''}{f_2'} ,
$$
so $f_1$ and $f_2$ differ by an affine transformation, implying $M_1 = M_2$.

\emph{Completeness}: Let $u \in C^{\infty}\left((0 , + \infty) , \R\right)$ be arbitrary. Let $v$ be an antiderivative of $u$, and let $w$ be an antiderivative of $e^{12 v}$. Then $w$ is strictly increasing and infinitely differentiable on $(0 , + \infty)$, and the integral mean of the first kind $M$ associated with $w$ satisfies (according to Proposition \ref{p11}):
$$
Q_M = \frac{1}{12} \dfrac{w''}{w'} = u .
$$
This proves completeness and achieves this proof.
\end{proof}

We now present some applications of the characteristic function of integral means of the first kind:

\begin{thm}\label{t4}
The only integral means of the first kind that are homogeneous are the Stolarsky means.
\end{thm}

\begin{proof}
We have already seen that the Stolarsky means are homogeneous and are also integral means of the first kind. Conversely, suppose $M$ is a homogeneous integral mean of the first kind. By Proposition \ref{p23}, the characteristic function of $M$ has the form $Q_M(x) = \frac{c}{x}$ for some $c \in \R$. On the other hand, if $f$ is an associated integrality function, Proposition \ref{p11} gives $Q_M = \frac{1}{12} \frac{f''}{f'}$. thus,
$$
\frac{1}{12} \dfrac{f''}{f'} = \dfrac{c}{x} ,
$$
which integrates to
$$
f(x) = \begin{cases}
\alpha x^{12 c + 1} + \beta & \text{if } c \neq - \frac{1}{12} \\
\alpha \log{x} + \beta & \text{if } c = - \frac{1}{12}
\end{cases} \quad (\alpha , \beta \in \R , \alpha \neq 0) .
$$
Replacing $f$ with an affine transformation of $f$ does not change $M$, so we may assume:
$$
f(x) = \begin{cases}
x^{12 c + 1} & \text{if } c \neq - \frac{1}{12} \\
\log{x} & \text{if } c = - \frac{1}{12}
\end{cases} .
$$
Thus, $M$ is the Stolarsky mean $S_{12 c + 2}$, completing this proof.
\end{proof}

\begin{coll}\label{c11}
The harmonic mean $H$ is not an integral mean of the first kind.
\end{coll}

\begin{proof}
Assume, for contradiction, that $H$ is an integral mean of the first kind. Being homogeneous, Theorem \ref{t4} implies that $H = S_p$ for some $p \in \R$. Equating characteristic functions, we find $p = - 4$, so that $H = S_{-4}$. This yields the identity:
$$
\dfrac{2 x y}{x + y} = \left(- \frac{1}{4} \, \dfrac{x^{-4} - y^{-4}}{x - y}\right)^{- 1/5}  \quad (\forall x , y > 0) .
$$
But as $x \to + \infty$, the left-hand side tends to $2 y$, while the right-hand side tends to $+ \infty$. Contradiction! Hence, $H$ is not an integral mean of the first kind.
\end{proof}

\begin{rmk}
We will see in the next subsection that $H$ belongs to the class of integral means of the second kind.
\end{rmk}

A further significant application of the characteristic function of integral means of the first kind lies in obtaining an approximation formula for definite integrals with sufficiently close bounds. We have the following proposition:

\begin{prop}\label{p12}
Let $f$ be a strictly monotonic and infinitely differentiable function on $(0 , + \infty)$. Then, for all sufficiently close positive real numbers $x$ and $y$, with $x < y$, we have
$$
\int_{x}^{y} f(t) \, d t \simeq (y - x) \, f\!\left(\dfrac{x \vabs{f'(x)}^{1/6} + y \vabs{f'(y)}^{1/6}}{\vabs{f'(x)}^{1/6} + \vabs{f'(y)}^{1/6}}\right)
$$
(with an error of order $O\left((x - y)^5\right)$).
\end{prop}

\begin{proof}
Let $M$ be the integral mean of the first kind having $f$ as an associated integrality function, and $N$ be the normal mean having ${\vabs{f'}}^{1/6}$ (which is infinitely differentiable, since $f'$ has a constant sign) as an associated weight function. By Propositions \ref{p6} and \ref{p11}, we have
$$
Q_M = Q_N = \frac{1}{12} \dfrac{f''}{f'} ,
$$
implying that $M$ and $N$ are close. More precisely, we have near the first bisector:
$$
M(x , y) = N(x , y) + O\left((x - y)^4\right) ;
$$
that is,
$$
f^{-1} \left(\frac{1}{y - x} \int_{x}^{y} f(t) \, d t\right) = \dfrac{x \vabs{f'(x)}^{1/6} + y \vabs{f'(y)}^{1/6}}{\vabs{f'(x)}^{1/6} + \vabs{f'(y)}^{1/6}} + O\left((x - y)^4\right) .
$$
Hence the required approximation formula of the proposition.
\end{proof}

\begin{rmk}
The problems of identifying the intersection of the class of integral means of the first kind with the class of normal means and with the class of additive means are still open.
\end{rmk}

\subsection{Integral means of the second kind}\label{subsec3-4}

To the best of our knowledge, this class has not been previously introduced in the literature. It arises from the idea of reversing the order of composition between $M$ and $f$ in the formula
(equivalent to \eqref{eq19})
$$
(f \circ M)(x , y) = \frac{1}{y - x} \int_{x}^{y} f(t) \, d t ,
$$
which defines an integral mean of the first kind.

\begin{defi}
A mean $M$ is called an \emph{integral mean of the second kind} if there exists a strictly monotonic and infinitely differentiable function $f$ on an interval $I \subseteq \R$, with $f(I) = (0 , + \infty)$, such that for all $x , y \in I$, with $x \neq y$, we have:
\begin{equation}\label{eq21}
M\left(f(x) , f(y)\right) = \frac{1}{y - x} \int_{x}^{y} f(t) \, d t .
\end{equation}
Such a function $f$ is called an \emph{integrality function} associated with $M$.
\end{defi}

\begin{rmk}\label{r4}
As in Remark \ref{r3}, performing the change of variable $t = x u + y (1 - u)$ in the integral on the right-hand side of \eqref{eq21} yields:
\begin{equation}\label{eq22}
M\left(f(x) , f(y)\right) = \int_{0}^{1} f\left(x u + y (1 - u)\right) \, d u .
\end{equation}
This formula also holds for $y = x$, demonstrating that any continuous function on $(0 , + \infty)^2$ of the form \eqref{eq21} (for $x \neq y$) defines a mean. Furthermore, Formula \eqref{eq22} facilitates the computation of the characteristic function for integral means of the second kind, and provides a basis for their extension to a broader class of means (see \S \ref{subsec3-6}).
\end{rmk}

One can readily verify that the arithmetic mean $A$, the geometric mean $G$, and the harmonic mean $H$ are all integral means of the second kind, with corresponding integrality functions given respectively by $f_1(t) = t$, $f_2(t) = \frac{1}{t^2}$, and $f_3(t) = \frac{1}{\sqrt{t}}$ ($t \in (0 , + \infty)$). It is also straightforward to see that if $f$ is an integrality function associated with an integral mean of the second kind $M$, then so is any function of the form $g(t) = f(\alpha t + \beta)$ ($\alpha , \beta \in \R$, $\alpha \neq 0$). The converse will be shown in Corollary \ref{c12} below.

The expression of the characteristic function of an integral mean of the second kind in terms of its associated integrality function is given by the following proposition:

\begin{prop}\label{p13}
Let $I$ be a nonempty open interval of $\R$ and $f : I \rightarrow (0 , + \infty)$ be a strictly monotonic and infinitely differentiable function, with $f(I) = (0 , + \infty)$. Let also $M$ be the integral mean of the second kind associated with $f$. Then, we have 
$$
Q_M = \frac{1}{6} \dfrac{(f^{-1})''}{(f^{-1})'} .
$$
\end{prop}

\begin{proof}
Differentiating Formula \eqref{eq22} twice with respect to $x$, we find (for all $x , y \in I$):
$$
f''(x) \dfrac{\partial M}{\partial x}\left(f(x) , f(y)\right) + f'(x)^2 \, \dfrac{\partial^2 M}{{\partial x}^2}\left(f(x) , f(y)\right) = \int_{0}^{1} u^2 f''\left(x u + y (1 - u)\right) \, d u .
$$
Then, setting $y = x$ and using Formula \eqref{eq-qay-1} gives
$$
Q_M\left(f(x)\right) = - \frac{1}{6} \dfrac{f''(x)}{f'^2(x)} \quad (\forall x \in I) .
$$
Substituting $x = f^{-1}(u)$ ($u \in (0 , + \infty)$), we find: 
$$
Q_M(u) = \left(- \frac{1}{6} \dfrac{f''}{f'^2} \circ f^{-1}\right)(u) = \frac{1}{6} \dfrac{(f^{-1})''}{(f^{-1})'}(u) ,
$$
as desired. 
\end{proof}

From Proposition \ref{p13}, we derive the following corollaries:

\begin{coll}\label{c12}
Let $f$ and $g$ be two integrality functions associated with the same integral mean of the second kind M. Then there exist $\alpha , \beta \in \R$, with $\alpha \neq 0$, such that $g(x) = f(\alpha x + \beta)$.
\end{coll}

\begin{proof}
Apply Proposition \ref{p13} to obtain:
$$
\dfrac{(f^{-1})''}{(f^{-1})'} = \dfrac{(g^{-1})''}{(g^{-1})'} ,
$$
implying that $\frac{(f^{-1})'}{(g^{-1})'}$ is constant (nonzero). Therefore, $f^{-1} = \alpha g^{-1} + \beta$ for some $(\alpha , \beta) \in \R^* \times \R$. Thus $f \circ f^{-1} \circ g = f \circ (\alpha g^{-1} + \beta) \circ g$; that is, $g(x) = f(\alpha x + \beta)$, as required.
\end{proof}

\begin{coll}\label{c13}
The class of integral means of the second kind is minimal and complete.
\end{coll}

\begin{proof}
\emph{Minimality}: Let $M_1$ and $M_2$ be two integral means of the second kind with the same characteristic function, and let $f_1$ and $f_2$ be their respective associated integrality functions. By Proposition \ref{p13}, we have
$$
\frac{1}{6} \dfrac{(f_1^{-1})''}{(f_1^{-1})'} = \frac{1}{6} \dfrac{(f_2^{-1})''}{(f_2^{-1})'} .
$$
As in the proof of Corollary \ref{c12}, it follows that there exists an affine function $\sigma$ such that $f_2 = f_1 \circ \sigma$. Therefore, $M_1 = M_2$.

\emph{Completeness}: Let $u \in C^{\infty}\left((0 , + \infty) , \R\right)$ be arbitrary. Let $v$ be an antiderivative of $u$ and $w$ be an antiderivative of $e^{6 v}$. Then, $w$ is strictly increasing and infinitely differentiable on $(0 , + \infty)$, implying that $w$ is a bijection from $(0 , + \infty)$ onto $I := w((0 , + \infty))$. Setting $f := w^{-1} : I \rightarrow (0 , + \infty)$, which is also strictly increasing and infinitely differentiable, the integral mean of the second kind $M$ associated with $f$ satisfies, by Proposition \ref{p13}:
$$
Q_M = \frac{1}{6} \dfrac{(f^{-1})''}{(f^{-1})'} = \frac{1}{6} \dfrac{w''}{w'} = u .
$$
This confirms completeness and concludes the proof.
\end{proof}

We conclude with the following proposition, which provides an approximation formula for definite integrals with sufficiently close bounds. Compared to its analogue, Proposition \ref{p12}, it is notably simpler.

\begin{prop}\label{p14}
Let $f$ be a strictly monotonic and infinitely differentiable function on an interval $I \subseteq \R$ such that $f(I) = (0 , + \infty)$. Then, for all sufficiently close real numbers $x$ and $y$ in $I$, with $x < y$, we have
$$
\int_{x}^{y} f(t) \, d t \simeq (y - x) \dfrac{f(x) \vabs{f'(y)}^{1/3} + f(y) \vabs{f'(x)}^{1/3}}{\vabs{f'(x)}^{1/3} + \vabs{f'(y)}^{1/3}}
$$
(with an error of order $O((x - y)^5)$).
\end{prop}

\begin{proof}
Let $M$ be the integral mean of the second kind having $f$ as an associated integrality function, and $N$ be the normal mean having ${\vabs{(f^{-1})'}}^{1/3}$ as an associated weight function. By Propositions \ref{p6} and \ref{p13}, we have
$$
Q_M = Q_N = \frac{1}{6} \dfrac{(f^{-1})''}{(f^{-1})'} ,
$$
implying that $M$ and $N$ are close. More precisely, we have near the first bisector:
$$
M(u , v) = N(u , v) + O\left((u - v)^4\right) .
$$
Setting $u = f(x)$ and $v = f(y)$, where $x , y \in I$, sufficiently close, we derive that:
$$
\frac{1}{y - x} \int_{x}^{y} f(t) \, d t = N\left(f(x) , f(y)\right) + O\left((x - y)^4\right) ,
$$
and the approximation follows from simplifying the expression of $N(f(x) , f(y))$.
\end{proof}

\begin{rmk}
The condition $f(I) = (0 , + \infty)$ in Proposition \ref{p14} can be omitted by considering means defined on $f(I)^2$ instead of $(0 , + \infty)^2$.
\end{rmk}

\begin{rmk}
The problem of describing the intersections between the class of integral mans of the second kind and other classes of means remains open.
\end{rmk}

\subsection{Weighted integral means of the first kind}\label{subsec3-5}

This class naturally extends the class of integral means of the first kind (see \S \ref{subsec3-3}). The idea is to consider, more generally than \eqref{eq20}, the function $M : (0 , + \infty)^2 \rightarrow (0 , + \infty)$ of the form:
$$
M(x , y) = f^{-1}\left(\int_{0}^{1} f\left(x t + y (1 - t)\right) g(t) \, d t\right) ,
$$
where $f$ is a strictly monotonic and infinitely differentiable function on $(0 , + \infty)$, and $g : [0 , 1] \rightarrow \R$ is a function satisfying minimal conditions to ensure that $M$ defines a mean. These conditions are precised in the following definition:

\begin{defi}\label{defi7}
A function $g : [0 , 1] \rightarrow \R$ is called a \emph{symmetric distribution function} if:
\begin{enumerate}
\item[(i)] $g$ is nonnegative almost everywhere on $[0 , 1]$.
\item[(ii)] $g$ is Lebesgue-integrable on $[0 , 1]$ and satisfies $\int_{0}^{1} g(t) \, d t = 1$.
\item[(iii)] $g(t) = g(1 - t)$ for almost all $t \in [0 , 1]$. 
\end{enumerate}
\end{defi}

\begin{nota}\label{nota1}
Given a strictly monotonic and infinitely differentiable function $f$ on $(0 , + \infty)$ and a symmetric distribution function $g$, we denote by $\I_{f , g}$ the real-valued function defined on $(0 , + \infty)^2$ by:
\begin{equation}\label{eq23}
\I_{f , g}(x , y) := f^{-1}\left(\int_{0}^{1} f\left(t x + (1 - t) y\right) g(t) \, d t\right) \quad (\forall x , y > 0) . 
\end{equation}
\end{nota}

We have the following proposition:

\begin{prop}\label{p15}
Let $f$ and $g$ be as in Notation \ref{nota1}. Then the function $\I_{f , g}$ defines a mean.
\end{prop}

\begin{proof}
Since $f$ is $C^{\infty}$ on $(0 , + \infty)$ then $\I_{f , g}$ is of class $C^{\infty}$ on $(0 , + \infty)^2$. Now, let $x , y \in (0 , + \infty)$ with $x \leq y$. Replacing $f$ by $(- f)$ if necessary (which leaves $\I_{f , g}$ unchanged), we may assume that $f$ is strictly increasing. Then, for all $t \in [0 , 1]$, we have:
$$
f(x) \leq f\left(t x + (1 - t) y\right) \leq f(y) .
$$
Multiplying by the function $g(t)$ and integrating over $[0 , 1]$ yields:
$$
f(x) \leq \int_{0}^{1} f\left(t x + (1 - t) y\right) g(t) \, d t \leq f(y) .
$$
Applying $f^{-1}$ gives $x \leq \I_{f , g}(x , y) \leq y$, as required. The symmetry of $\I_{f , g}$ follows by changing variables $u = 1 - t$ in \eqref{eq23} and using $g(1 - u) = g(u)$ for almost all $u \in [0 , 1]$. This completes the proof.
\end{proof}

\noindent\textbf{Terminology:} In the context of Proposition \ref{p15}, the mean $\I_{f , g}$ is called the \emph{weighted integral mean of the first kind} with \emph{integrality function} $f$ and \emph{symmetric distribution function} $g$. When $g$ is fixed, a mean of the from $\I_{f , g}$ is referred to as the \emph{$g$-weighted integral mean of the first kind} with \emph{integrality function} $f$.

\begin{rmks}\label{rmks1} The following properties are straightforward to verify:
\begin{enumerate}
\item The $1$-weighted integral means of the first kind are precisely the integral means of the first kind discussed in Subsection \ref{subsec3-3}.
\item For any symmetric distribution function $g$, we have:
$$
\int_{0}^{1} t g(t) \, d t = \int_{0}^{1} (1 - t) g(t) \, d t = \frac{1}{2} ,
$$
implying in particular that $\I_{\id_{(0 , + \infty)} , g} = A$.
\item For every symmetric distribution function $g$ and every odd positive integer $n$, we have:
$$
\int_{0}^{1} \left(t - \frac{1}{2}\right)^n g(t) \, d t = 0 .
$$
\item Let $g$ be a symmetric distribution function. If $f$ is an integrality function associated with a $g$-weighted integral mean of the first kind, then so is any affine transformation $\alpha f + \beta$ ($\alpha , \beta \in \R$, $\alpha \neq 0$). The converse of this property also holds and can be readily proved using Proposition \ref{p17} below.
\end{enumerate}
\end{rmks}

\subsubsection*{Nontrivial examples of symmetric distribution functions}
\addcontentsline{toc}{subsubsection}{Nontrivial examples of symmetric distribution functions}

A nontrivial example of a symmetric distribution function can be constructed by selecting a function $h : [0 , 1/4] \rightarrow \R$ that satisfies the following conditions:
\begin{itemize}
\item $h$ is nonnegative almost everywhere on $[0 , 1/4]$.
\item The function $u \mapsto \frac{h(u)}{\sqrt{1 - 4 u}}$ is Lebesgue-integrable on $[0 , 1/4]$ with a nonzero integral.
\end{itemize}
Then, the function $g : [0 , 1] \rightarrow \R$, defined by
$$
g(t) := \dfrac{h\left(t (1 - t)\right)}{\int_{0}^{1} h\left(x (1 - x)\right) \, d x} \quad (\forall t \in [0 , 1]) ,
$$
is a symmetric distribution function. Notably,
$$
\int_{0}^{1} h\left(x (1 - x)\right) \, d x = 2 \int_{0}^{1/2} h\left(x (1 - x)\right) \, d x = 2 \int_{0}^{1/4} \frac{h(u)}{\sqrt{1 - 4 u}} \, d u .
$$
Choosing $h(t) = t^r$ with $r > - 1$ yields the symmetric distribution function $\mathfrak{g}_r$, given by:
$$
\mathfrak{g}_r(t) = \dfrac{t^r (1 - t)^r}{\beta(r + 1 , r + 1)} \quad (\forall t \in [0 , 1]) ,
$$
where $\beta$ denotes the Euler Beta function. The particular case $r = - \frac{1}{2}$ gives the following important symmetric distribution function:
\begin{equation}\label{eq24}
\mathfrak{g}(t) = \dfrac{1}{\pi \sqrt{t (1 - t)}} \quad (\forall t \in (0 , 1)) .
\end{equation}
We will see later that the Gauss $\AGM$ mean is closely related to the class of $\mathfrak{g}$-weighted integral means of the first kind (see \S \ref{sec4}).

\subsubsection*{Essential integrals of a symmetric distribution function}
\addcontentsline{toc}{subsubsection}{Essential integrals of a symmetric distribution function}

Given a symmetric distribution function $g$, we will see below that the dependence on $g$ of the characteristic functions (of different orders) of a $g$-weighted integral mean of the first kind is expressed solely in terms of certain integrals related to $g$, which we have termed the \emph{essential integrals} of $g$.

\begin{defi}\label{defi8}
Let $g$ be a symmetric distribution function and $n$ be a nonnegative integer. We define the \emph{essential integral} of $g$ of order $n$ as the positive real number $c_n(g)$ given by:
$$
c_n(g) := \int_{0}^{1} \left(t - \frac{1}{2}\right)^{\!\! 2 n} g(t) \, d t .
$$
In particular, $c_0(g) = 1$.
\end{defi}

We have the following proposition:

\begin{prop}\label{p16}
Let $f$ be a function expandable in a Taylor series around every point of $(0 , + \infty)$, and let $g$ be a symmetric distribution function. Then, near the first bisector, we have:
$$
\int_{0}^{1} f\left(x t + y (1 - t)\right) g(t) \, d t = f(A) + \frac{f''(A)}{2!} c_1(g) (x - y)^2 + \frac{f^{(4)}(A)}{4!} c_2(g) (x - y)^4 + \dots .
$$
\end{prop}

\begin{proof}
Let $x , y > 0$ be sufficiently close. The Taylor expansion of the real-valued function $t \rightarrow f(x t + y (1 - t))$ around $t = 1/2$ is:
$$
f\left(x t + y (1 - t) y\right) = f(A) + (x - y) \frac{f'(A)}{1!} \left(t - \frac{1}{2}\right) + (x - y)^2 \frac{f''(A)}{2!} \left(t - \frac{1}{2}\right)^2 + \dots .
$$
Multiplying by $g(t)$ and integrating over $[0 , 1]$ yields
$$
\int_{0}^{1} f\left(x t + y (1 - t)\right) g(t) \, d t = \sum_{n = 0}^{+ \infty} \frac{f^{(n)}(A)}{n!} (x - y)^n \int_{0}^{1} \left(t - \frac{1}{2}\right)^n g(t) \, d t .
$$
Since $\int_{0}^{1} (t - 1/2)^n g(t) \, d t = 0$ for odd $n$ (see Item 3 of Remarks \ref{rmks1}), we conclude that
$$
\int_{0}^{1} f\left(x t + y (1 - t)\right) g(t) \, d t = \sum_{n = 0}^{+ \infty} \frac{f^{(2 n)}(A)}{(2 n)!} (x - y)^{2 n} \int_{0}^{1} \left(t - \frac{1}{2}\right)^{2 n} g(t) \, d t ,
$$
as required.
\end{proof}

\begin{rmk}
Proposition \ref{p16} shows that if $f$ is strictly monotonic on $(0 , + \infty)$ and has a Taylor expansion around every point of $(0 , + \infty)$, and $g$ is a symmetric distribution function, then the expansion of the mean $\I_{f , g}$ (near the first bisector) depends only on $f$ and on the essential integrals of $g$. 
\end{rmk}

In the following proposition, we explicitly express the characteristic function of a weighted integral mean of the first kind in terms of its integrality function and its symmetric distribution function.

\begin{prop}\label{p17}
Let $f$ be a strictly monotonic and infinitely differentiable function on $(0 , + \infty)$, and $g : [0 , 1] \rightarrow \R$ be a symmetric distribution function. Then the characteristic function of the $g$-weighted integral mean of the first kind $\I_{f , g}$ is given by:
$$
Q_{\I_{f , g}} = c_1(g) \dfrac{f''}{f'} .
$$
\end{prop} 

\begin{proof}
Set $M := \I_{f , g}$. From Formula \eqref{eq23}, we have for all $x , y > 0$:
$$
f\left(M(x , y)\right) = \int_{0}^{1} f\left(x t + y (1 - t)\right) g(t) \, d t .
$$
Differentiating twice with respect to $x$ yields
$$
\dfrac{\partial^2 M}{{\partial x}^2}(x , y) f'\left(M(x , y)\right) + \left(\dfrac{\partial M}{\partial x}(x , y)\right)^2 f''\left(M(x , y)\right) = \int_{0}^{1} t^2 f''\left(x t + y (1 - t)\right) g(t) \, d t .
$$
By setting $y = x$ and using Formula \eqref{eq-qay-1}, we derive that
$$
Q_M(x) = \left(\int_{0}^{1} t^2 g(t) \, d t - \frac{1}{4}\right) \dfrac{f''(x)}{f'(x)} .
$$
Finally, noting that
\begin{multline*}
\int_{0}^{1} t^2 g(t) \, d t = \int_{0}^{1} \left[\frac{1}{4} + \left(t - \frac{1}{2}\right) + \left(t - \frac{1}{2}\right)^2\right] g(t) \, d t \\
= \frac{1}{4} \int_{0}^{1} g(t) \, d t + \int_{0}^{1} \left(t - \frac{1}{2}\right) g(t) \, d t + \int_{0}^{1} \left(t - \frac{1}{2}\right)^2 g(t) \, d t = c_1(g) + \frac{1}{4}  
\end{multline*}
(in view of Item 3 of Remarks \ref{rmks1}), the required result follows.
\end{proof}

\begin{rmk}
Taking $g(t) = 1$ in Proposition \ref{p17} recovers Proposition \ref{p11}.  
\end{rmk}

From Proposition \ref{p17}, we readily derive the following result.

\begin{coll}
For every fixed symmetric distribution function $g$, the class of $g$-weighted integral means of the first kind is minimal and complete.
\end{coll}

\begin{proof}
Similar to that of Corollary \ref{c10}.
\end{proof}

\subsection{Weighted integral means of the second kind}\label{subsec3-6}

This class naturally extends the class of integral means of the second kind by weighting the integral in its definition with a symmetric distribution function. We begin with the following proposition:

\begin{prop}\label{p18}
Let $f$ be a strictly monotonic and infinitely differentiable function on an interval $I \subseteq \R$, with $f(I) = (0 , + \infty)$, and let $g$ be a symmetric distribution function. Then there exists a unique mean $M$ satisfying, for all $x , y \in I$:
\begin{equation}\label{eq25}
M\left(f(x) , f(y)\right) = \int_{0}^{1} f\left(x t + y (1 - t)\right) g(t) \, d t .
\end{equation} 
\end{prop}

\begin{proof}
Similar to that of Proposition \ref{p15}.
\end{proof}

\begin{nota}\label{nota2}
In the context of Proposition \ref{p18}, we denote by $\I_{f , g}^*$ the unique mean defined through Formula \eqref{eq25}. More explicitly, for all $x , y > 0$, we have:
\begin{equation}
\I_{f , g}^*(x , y) := \int_{0}^{1} f\left(t f^{-1}(x) + (1 - t) f^{-1}(y)\right) g(t) \, d t .
\end{equation}
\end{nota}

\noindent\textbf{Terminology:} In the context of Proposition \ref{p18} and Notation \ref{nota2}, the mean $\I_{f , g}^*$ is called the \emph{weighted integral mean of the second kind} with \emph{integrality function} $f$ and \emph{symmetric distribution function} $g$. When $g$ is fixed, a mean of the form $\I_{f , g}^*$ is referred to as the \emph{$g$-weighted integral mean of the second kind} with \emph{integrality function} $f$. 

\begin{rmks}\label{rmks2}
The following properties are straightforward to verify:
\begin{enumerate}
\item The $1$-weighted integral means of the second kind are precisely the integral means of the second kind discussed in Subsection \ref{subsec3-4}.
\item For any symmetric distribution function $g$, we have:
$$
\I_{\id_{(0 , + \infty)} , g}^* = A .
$$
\item Let $g$ be a symmetric distribution function. If $f$ is an integrality function associated with a $g$-weighted integral mean of the second kind, then so is any function of the form $\widetilde{f}(t) = f(\alpha t + \beta)$ ($\alpha , \beta \in \R$, $\alpha \neq 0$). The converse of this property also holds and can be readily proved using Proposition \ref{p19} below.
\end{enumerate}
\end{rmks}

We now proceed to express the characteristic function of a weighted integral mean of the second kind in terms of its integrality function and its symmetric distribution function. We have the following proposition:

\begin{prop}\label{p19}
Let $f$ and $g$ be as in Proposition \ref{p18}. Then the characteristic function of the $g$-weighted integral mean of the second kind $\I_{f , g}^*$ is given by:
$$
Q_{\I_{f , g}^*} = \mu(g) \dfrac{(f^{-1})''}{(f^{-1})'} ,
$$
where
$$
\mu(g) := \int_{0}^{1} t (1 - t) g(t) \, d t > 0 .
$$
\end{prop}

\begin{proof}
Consider the context of Proposition \ref{p18}, where $M = \I_{f , g}^*$. Differentiating Formula \eqref{eq25} twice with respect to $x$, we obtain (for all $x , y \in I$):
$$
f''(x) \dfrac{\partial M}{\partial x}\left(f(x) , f(y)\right) + f'(x)^2 \dfrac{\partial^2 M}{{\partial x}^2}\left(f(x) , f(y)\right) = \int_{0}^{1} t^2 f''\left(x t + y (1 -  t)\right) g(t) \, d t .
$$
Setting $y = x$ and using Formula \eqref{eq-qay-1} yields
$$
Q_M\left(f(x)\right) = \left(\int_{0}^{1} t^2 g(t) \, d t - \frac{1}{2}\right) \dfrac{f''(x)}{f'^2(x)} .
$$
Substituting $x = f^{-1}(u)$ ($u \in (0 , + \infty)$) and noting that
$$
\int_{0}^{1} t^2 g(t) \, d t - \frac{1}{2} = \int_{0}^{1} t (t - 1) g(t) \, d t , 
$$
we get
$$
Q_M(u) = \left(\int_{0}^{1} t (t - 1) g(t) \, d t\right) \left(\dfrac{f''}{f'^2} \circ f^{-1}\right)(u) = \left(\int_{0}^{1} t (1 - t) g(t) \, d t\right) \dfrac{(f^{-1})''}{(f^{-1})'}(u) ,
$$
as required.
\end{proof}

From Proposition \ref{p19}, we derive the following corollary:

\begin{coll}\label{c14}
For every fixed symmetric distribution function $g$, the class of $g$-weighted integral means of the second kind is minimal and complete.
\end{coll}

\begin{proof}
Similar to that of Corollary \ref{c13}.
\end{proof}

\section{$M$-means}\label{sec4}

Given a mean $M$, we can construct related means by composing an appropriate regular function $f$ on the right and its inverse $f^{-1}$ on the left. This leads to the concept of \emph{$M$-means}.

\begin{defi}\label{defi9}
Let $M$ be a mean. An \emph{$M$-mean} is any function $F : (0 , + \infty)^2 \rightarrow (0 , + \infty)$ of the form
$$
F = f^{-1} \circ M \circ f ,
$$
where $f : (0 , + \infty) \rightarrow (0 , + \infty)$ is a strictly monotonic and infinitely differentiable function. Clearly, $F$ defines a mean.
\end{defi}

\begin{nota}
Given a mean $M$, we denote by $\mathscr{C}_M$ the class of all $M$-means.
\end{nota}

\begin{rmk}
The class $\mathscr{C}_A$ of $A$-means si strictly contained in the class of additive means discussed in \S \ref{subsec3-2}. For instance, the geometric mean $G$ is additive but not an $A$-mean, because in the expression $(\log)^{-1} \circ A \circ (\log)$ of $G$, the logarithmic function $\log$ takes nonpositive values, violating the domain requirements of $f$ in Definition \ref{defi9}.
\end{rmk}

A central question is to determine, for a given mean $M$, when the class $\mathscr{C}_M$ is complete and when it is minimal. We partially answer this question below in the most frequent cases, after first establishing a closed-from expression for the characteristic function of $M$-means. We have the following proposition:

\begin{prop}\label{p20}
Let $M$ be a mean, and let $f : (0 , + \infty) \rightarrow (0 , + \infty)$ be a strictly monotonic and infinitely differentiable function. Then the characteristic function of the $M$-mean $(f^{-1} \circ M \circ f)$ is given by:
$$
Q_{f^{-1} \circ M \circ f} = \frac{1}{4} \dfrac{f''}{f'} + f' \cdot \left(Q_M \circ f\right) .
$$
\end{prop}

\begin{proof}
Setting $\widetilde{M} := f^{-1} \circ M \circ f$, we have for all $x , y > 0$:
$$
f\left(\widetilde{M}(x , y)\right) = M\left(f(x) , f(y)\right) .
$$
Differentiating twice with respect to $x$, we get (for all $x , y > 0$):
\begin{multline*}
\dfrac{\partial^2 \widetilde{M}}{{\partial x}^2}(x , y) \cdot f'\left(\widetilde{M}(x , y)\right) + \left(\dfrac{\partial \widetilde{M}}{\partial x}(x , y)\right)^2 f''\left(\widetilde{M}(x , y)\right) = f''(x) \cdot \dfrac{\partial M}{\partial x}\left(f(x) , f(y)\right) \\
+ f'(x)^2 \cdot \dfrac{\partial^2 M}{{\partial x}^2}\left(f(x) , f(y)\right) .
\end{multline*}
Setting $y = x$ and using Formula \eqref{eq-qay-1} yields the required result.
\end{proof}

For what follows, we adopt the notation $\int u$ to denote an antiderivative of a given continuous real-valued function $u$ on $(0 , + \infty)$.

From Proposition \ref{p20}, we derive the following corollary:

\begin{coll}\label{c15}
Let $M$ be a mean such that, for some choice of antiderivatives, we have
$$
\im\left(\int e^{4 \int Q_M}\right) = (0 , + \infty) .
$$
Then, for all $u \in C^{\infty}\left((0 , + \infty) , \R\right)$ satisfying $\im\left(\int e^{4 \int u}\right) = (0 , + \infty)$ (for some choice of antiderivatives) there exists an $M$-mean $M_0$ such that $Q_{M_0} = u$. Furthermore, the class $\mathscr{C}_M$ is minimal if and only if $\mathscr{C}_M = \mathscr{C}_A$.
\end{coll}

The proof of Corollary \ref{c15} relies on the following two lemmas:

\begin{lemma}\label{l1}
Let $M$ be a mean such that, for some choice of antiderivatives, we have \linebreak $\im\left(\int e^{4 \int Q_M}\right) = (0 , + \infty)$. Then the $M$-mean $N$ given by:
$$
N := \left(\int e^{4 \int Q_M}\right) \circ M \circ \left(\int e^{4 \int Q_M}\right)^{-1} 
$$
is almost arithmetic (i.e., $Q_N = 0$).
\end{lemma}

\begin{proof}
Apply Proposition \ref{p20} for $f = \left(\int e^{4 \int Q_M}\right)^{-1}$.
\end{proof}

\begin{lemma}\label{l2}
Let $M$ be a mean. For all $\alpha > 0$ and all $\beta \geq 0$, let $L_{\alpha , \beta}$ denote the affine function defined on $(0 , + \infty)$ by:
$$
L_{\alpha , \beta}(x) := \alpha x + \beta \quad (\forall x > 0) .
$$
Suppose that we have for all $\alpha > 0$ and all $\beta \geq 0$:
$$
M \circ L_{\alpha , \beta} = L_{\alpha , \beta} \circ M .
$$
Then $M = A$.
\end{lemma}

\begin{proof}
The hypothesis for $\beta = 0$ means that $M$ is homogeneous. For $\alpha = 1$, the hypothesis means that $M$ is additively homogeneous. Therefore, $M$ is both homogeneous and additively homogeneous, so by Proposition \ref{p1}, $M = A$, as required. 
\end{proof}

\begin{proof}[Proof of Corollary \ref{c15}]
If necessary, replace $M$ with the $M$-mean
$$
N := \left(\int e^{4 \int Q_M}\right) \circ M \circ \left(\int e^{4 \int Q_M}\right)^{-1} ,
$$
which is almost arithmetic by Lemma \ref{l1} and satisfies $\mathscr{C}_N = \mathscr{C}_M$. So we may assume that $M$ is itself almost arithmetic; i.e., $Q_M = 0$. Let $u \in C^{\infty}\left((0 , + \infty) , \R\right)$ satisfying $\im\left(\int e^{4 \int u}\right) = (0 , + \infty)$ for some choice of antiderivatives, and define $f := \left(\int e^{4 \int u}\right)$. Then, $f$ is positive, strictly increasing, and infinitely differentiable on $(0 , + \infty)$. Moreover, $\frac{1}{4} \frac{f''}{f'} = u$. It follows from Proposition \ref{p20} that:
\begin{align*}
Q_{f^{-1} \circ M \circ f} & = \frac{1}{4} \dfrac{f''}{f'} + f' \cdot \left(Q_M \circ f\right) \\
& = u \quad (\text{since } Q_M = 0 \text{ and } \frac{1}{4} \frac{f''}{f'} = u) .
\end{align*}
Hence, the $M$-mean $M_0 := f^{-1} \circ M \circ f$ satisfies $Q_{M_0} = u$. 

Now, let us prove the second part of Corollary \ref{c15} concerning the minimality of the class $\mathscr{C}_M$. Since $\mathscr{C}_A$ is a subclass of the class of additive means, which is minimal by Corollary \ref{c7}, then $\mathscr{C}_A$ is minimal. Conversely, suppose that $\mathscr{C}_M \neq \mathscr{C}_A$ and show that $\mathscr{C}_M$ is not minimal. Since $M \neq A$, then by Lemma \ref{l2}, there exist $\alpha > 0$ and $\beta \geq 0$ such that for the affine function $L(x) = \alpha x + \beta$, we have
$$
M \circ L \neq L \circ M .
$$
Consequently, the $M$-mean $M^* := L^{-1} \circ M \circ L$ is different from $M$, while according to Proposition \ref{p20}:
\begin{align*}
Q_{M^*} & = \frac{1}{4} \dfrac{L''}{L'} + L' \cdot \left(Q_M \circ L\right) \\
& = 0 \quad (\text{since } Q_M = 0 \text{ and } L'' = 0) \\
& = Q_M .
\end{align*}
This confirms that the class $\mathscr{C}_M$ is not minimal, completing the proof of the second part of Corollary \ref{c15} and thus the entire proof. 
\end{proof}

We conclude this section, and subsequently this paper, by relating the arithmetic-geometric mean $\AGM$ to the weighted integral means of the first kind via the concept of $M$-means. We have the following proposition:

\begin{prop}\label{p21}
Let $M$ be the weighted integral mean of the first kind with integrality function $t \mapsto \frac{1}{\sqrt{t}}$ and symmetric distribution function $t \mapsto \frac{1}{\pi \sqrt{t (1 - t)}}$. Then the $\AGM$ mean is an $M$-mean. More precisely, setting $h(t) = t^2$, we have:
$$
\AGM = h^{-1} \circ \I_{\frac{1}{\sqrt{t}} , \frac{1}{\pi \sqrt{t (1 - t)}}} \circ h .
$$
\end{prop}

\begin{proof}
As established in \S \ref{subsec3-5} (see \eqref{eq24}), the function $\mathfrak{g} : (0 , 1) \rightarrow \R$ defined by 
$$
\mathfrak{g}(t) := \frac{1}{\pi \sqrt{t (1 - t)}} \quad (\forall t \in (0 , 1))
$$
is a symmetric distribution function. Now, given $x , y > 0$, we compute:
\begin{align*}
\left(h^{-1} \circ \I_{\frac{1}{\sqrt{t}} , \mathfrak{g}(t)} \circ h\right)(x , y) & = \sqrt{\I_{\frac{1}{\sqrt{t}} , \mathfrak{g}(t)}\left(x^2 , y^2\right)} \\
& = \left(\int_{0}^{1} \dfrac{d t}{\sqrt{t x^2 + (1 - t) y^2} \cdot \pi \sqrt{t (1 - t)}}\right)^{-1} .
\end{align*}
Performing the change of variable $t = \cos^2\theta$ (with $\theta \in [0 , \frac{\pi}{2}]$) in the last integral yields:
$$
\left(h^{-1} \circ \I_{\frac{1}{\sqrt{t}} , \mathfrak{g}(t)} \circ h\right)(x , y) = \dfrac{\pi}{2} \left(\int_{0}^{\pi/2} \dfrac{d \theta}{\sqrt{x^2 \cos^2\theta + y^2 \sin^2\theta}}\right)^{-1} = \AGM(x , y)
$$
(according to Formula \eqref{eq7}). The proposition is proved.
\end{proof}

\begin{rmk}
The expression for the characteristic function of the $\AGM$ mean (already established in Corollary \ref{c5}) can be recovered by applying the formula from Proposition \ref{p21} in conjunction with those from Propositions \ref{p20} and \ref{p17}.
\end{rmk}

\rhead


{\it References}

\begin{thebibliography}{99}
\bibitem{alz}
{\sc A. Alzer \& S. Ruscheweyh}. On the intersection of two-parameter mean value families, {\it Proc. Amer. Math. Soc.}, {\bf 29} n°9 (2001), p. 2655-2662.

\bibitem{bor}
{\sc J. M. Borwein \& P. B. Borwein}. {\it Pi and the AGM}, A Study in Analytic Number Theory and Computational Complexity, Wiley, New York, 1987.

\bibitem{bul}
{\sc P. S. Bullen}. {\it Handbook of means and their inequalities}, {\bf 560}, Springer Science \& Business Media, 2013.

\bibitem{bu-el}
{\sc T. Buri\'{c} \& N. Elezovi\'{c}}. Computation and analysis of the asymptotic expansions of the compound means, {\it Appl. Math. Comput.}, {\bf 303} (2017), p. 48-54.

\bibitem{el1}
{\sc N. Elezovi\'{c}}. Asymptotic inequalities and comparison of classical means, {\it J. Math. Inequal.}, {\bf 9}, n°1 (2015), p. 177-196.

\bibitem{elv1}
{\sc N. Elezovi\'{c} \& L. Vuk\v{s}i\'{c}}. Asymptotic expansions of bivariate classical means and related inequalities, {\it J. Math. Inequal.}, {\bf 8}, n°4 (2014), p. 707-7024.

\bibitem{elv2}
{\sc N. Elezovi\'{c} \& L. Vuk\v{s}i\'{c}}. Asymptotic expansions and comparison of bivariate parameter means, {\it Math. Inequal. Appl.}, {\bf 17}, n°4 (2014), p. 1225-1244.

\bibitem{elv3}
{\sc N. Elezovi\'{c} \& L. Vuk\v{s}i\'{c}}. Asymptotic expansions of integral means and applications to the ratio of gamma functions, {\it Appl. Math. Comput.}, {\bf 235} (2014), p. 187-200.

\bibitem{far}
{\sc B. Farhi}. Algebraic and topological structures on the set of mean functions and generalization of the AGM mean, {\it Colloq. Math.}, {\bf 132} (2013), p. 139-149.

\bibitem{gini}
{\sc C. Gini}. Di una Formula Compressiva delle Medie, {\it Metron}, {\bf 13} (1938), p. 3-22.

\bibitem{leh}
{\sc D. H. Lehmer}. On the compounding of certain means, {\it J. Math. Anal. Appl.}, {\bf 36} (1971), p. 183-200.

\bibitem{mih}
{\sc L. Mihokovi\`c}. Coinciding mean of the two symmetries on the set of mean functions, {\it Axioms}, {\bf 12} n°3 (2023), 238. 

\bibitem{sto}
{\sc K. B. Stolarsky}. Generalizations of the logarithmic mean, {\it Math. Mag.}, {\bf 48}, n°2 (1975), p. 87-92.

\end{thebibliography}
\end{document}